\documentclass[11pt]{article}
\newcommand{\documentdate}{3 V 2022}

\usepackage{a4wide,latexsym,amsmath,varioref,graphicx,xcolor}
\usepackage{xcolor}
\title{Convergence properties of an Objective-Function-Free\\
       Optimization regularization algorithm, including \\
       an $\mathcal{O}(\epsilon^{-3/2})$ complexity bound}

\author{	
   S. Gratton%
   \thanks{Universit\'e de Toulouse, INP, IRIT, Toulouse, France. Email:
     serge.gratton@enseeiht.fr. Work partially supported by 3IA Artificial and
     Natural Intelligence Toulouse Institute (ANITI), French "Investing for the Future
     - PIA3" program under the Grant agreement ANR-19-PI3A-0004"}, 
   ~S. Jerad%
   \thanks{ANITI, Universit\'e de Toulouse, INP, IRIT, Toulouse, France. Email:
     sadok.jerad@toulouse-inp.fr}
   ~and Ph. L. Toint%
   \thanks{NAXYS, University of Namur, Namur, Belgium. Email:
     philippe.toint@unamur.be. Partly supported by ANITI.}
}

\newcommand{\beqn}[1]{\begin{equation}\label{#1}}
\newcommand{\eeqn}{\end{equation}}
\newcommand{\req}[1]{(\ref{#1})}
\newcommand{\ms}{\;\;\;\;}
\newcommand{\tim}[1]{\;\; \mbox{#1} \;\;}

\newtheorem{theorem}{Theorem}[section]
\newtheorem{lemma}[theorem]{Lemma}
\newtheorem{corollary}{Corollary}

\newcommand{\numsection}[1]{\section{#1}\setcounter{equation}{0}}
\newcommand{\appnumsection}[1]{\section*{#1}\setcounter{equation}{0}
  \renewcommand{\theequation}{A.\arabic{equation}}
  \renewcommand{\thetheorem}{A.\arabic{theorem}}
  \renewcommand{\thetable}{A.\arabic{table}}
  \renewcommand{\thefigure}{A.\arabic{figure}}
  \renewcommand{\thesection}{A} }\renewcommand{\theequation}{\arabic{section}.\arabic{equation}}

\newcounter{algo}[section]
\renewcommand{\thealgo}{\thesection.\arabic{algo}}
\newcommand{\llem}[2]{\vspace{\baselineskip} 
\noindent\framebox[\textwidth]{\parbox{0.95\textwidth}{
\begin{lemma} \label{#1} \rm #2 \end{lemma} } } \vspace{\baselineskip} }
\newcommand{\algo}[3]{\refstepcounter{algo}
\begin{center}\begin{figure}[htbp]
\framebox[\textwidth]{
\parbox{0.95\textwidth} {\vspace{\topsep}
{\bf Algorithm \thealgo : #2}\label{#1}\\
\vspace*{-\topsep} \mbox{ }\\
{#3} \vspace{\topsep} }}
\end{figure}\end{center}}
\newcommand{\bpr}{{\bf Proof.} \hspace{1.5mm}}
\newcommand{\epr}{\hfill $\Box$ \vspace*{1em}}
\newcommand{\proof}[1]{
\begin{list}{}{
\setlength{\topsep}{0.0pt}
\setlength{\partopsep}{0.0pt}
\setlength{\leftmargin}{0.025\textwidth}
\setlength{\rightmargin}{0.5\leftmargin}
\setlength{\labelwidth}{0.5\leftmargin}
\setlength{\labelsep}{0.25\leftmargin}}
\item \bpr #1 \epr \noindent
\end{list}}
\newcommand{\lthm}[2]{\vspace{\baselineskip} 
\noindent\framebox[\textwidth]{\parbox{0.95\textwidth}{
\begin{theorem} \label{#1} \rm #2 \end{theorem} } } \vspace{\baselineskip} }
\newcommand{\lcor}[2]{\vspace{\baselineskip} 
\noindent\framebox[\textwidth]{\parbox{0.95\textwidth}{
\begin{corollary} \label{#1} \rm #2 \end{corollary} } } \vspace{\baselineskip}
}

\newcommand{\iiz}[1]{\{ 0, \ldots, #1 \}}
\newcommand{\iibe}[2]{\{ #1, \ldots, #2 \}}

\newcommand{\calO}{{\cal O}}

\renewcommand{\Re}{\hbox{I\hskip -2pt R}}
\newcommand{\bigfrac}[2]{\frac{\displaystyle #1}{\displaystyle #2}}
\newcommand{\sfrac}[2]{{\scriptstyle \frac{#1}{#2}}}
\newcommand{\half}{\sfrac{1}{2}}
\newcommand{\eqdef}{\stackrel{\rm def}{=}}

\newcommand{\kap}[1]{\kappa_{\mbox{\tiny #1}}}
\newcommand{\khigh}{\kappa_{\mbox{\tiny high}}}
\newcommand{\al}[1]{{\footnotesize{\sf #1}}}
\newcommand{\tal}[1]{{\normalsize {\sf #1}}}
\newcommand{\ttal}[1]{{\large {\sf #1}}}

\newcommand{\comment}[1]{}
\newcommand{\private}[1]{}

\date{\documentdate}

\begin{document}

\maketitle

\begin{abstract}
An adaptive regularization algorithm for unconstrained nonconvex
optimization is presented in which the objective function is never
evaluated, but only derivatives are used. This algorithm belongs to
the class of adaptive regularization methods, for which optimal
worst-case complexity results are known for the standard framework
where the objective function is evaluated.  It is shown in this paper
that these excellent complexity bounds are also valid for the new
algorithm, despite the fact that significantly less
information is used. In particular, it is shown that, if derivatives of
degree one to $p$ are used, the algorithm will find an $\epsilon_1$-approximate
first-order minimizer in at most $\calO(\epsilon_1^{-(p+1)/p})$
iterations, and an $(\epsilon_1,\epsilon_2)$-approximate second-order
minimizer in at most $\calO(\max[\epsilon_1^{-(p+1)/p},\epsilon_2^{-(p+1)/(p-1)}])$
iterations. As a special case, the new algorithm using first
and second derivatives, when applied to functions with Lipschitz
continuous Hessian, will find an iterate $x_k$ at which the
gradient's norm is less than $\epsilon_1$ in at most
$\calO(\epsilon_1^{-3/2})$ iterations.
\end{abstract}

{\small
\textbf{Keywords:} nonlinear optimization, adaptive regularization methods,
evaluation complexity, \\objective-function-free optimization (OFFO).
}

\numsection{Introduction}

This paper is about the (complexity-wise) fastest known optimization
method which does not evaluate the objective function.  Such methods,
coined OFFO for Objective-Function-Free Optimization, have recently
been very popular in the context of noisy problems, in particular in
deep learning applications (see
\cite{KingBa15,DuchHazaSing11,WuWardBott18,WardWuBott19} among many others), where they have shown remarkable
insensitivity to the noise level. This is a first motivation to
consider them, and it is our point of view that their deterministic
(noiseless) counterparts are good stepping stones to understand their
behaviour. Another motivation is the observation that other more
standard methods (using objective function evaluations) have been
proposed in the noisy case, but these typically require the noise on the
function values to be tightly controlled at a level lower than
that allowed for derivatives
\cite{CartSche17,ChenMeniSche18,BlanCartMeniSche19,
  CurtScheShi19,BeraCaoSche21,BellGuriMoriToin21c,BellGuriMoriToin21,BellGuriMoriToin20b}

The convergence anaysis of OFFO algorithms is not a new subject, and
has been considered for instance in
\cite{DefoBottBachUsun20,WardWuBott19,GratJeraToin22a,GratJeraToin22b,
GrapStel22,WuWardBott18}.  However, as far as the authors
are aware, the existing theory focuses on the case where only
gradients are used (with the exception of \cite{GratToin22}) and
establish a worst-case iteration complexity of, at best,
$\calO(\epsilon^{-2})$ for finding an $\epsilon$-approximate
first-order stationary point \cite{Nest04}. It is already remarkable
that this bound is, in order and for the same goal, identical to that
of standard methods using function values. But methods using
second-derivatives have proved to be globally more efficient in this
latter context, and the (complexity-wise) fastest such method is known
to have an $\calO(\epsilon^{-3/2})$ complexity bound
\cite{NestPoly06,CurtRobiSama17,RoyeWrig18,CartGoulToin11d,BirgGardMartSantToin17,
CartGoulToin22}.
Moreover, this better bound was shown to be sharp and
optimal among a large class of optimization algorithms using
second-derivatives for the noiseless case \cite{CartGoulToin18a}.
Is such an improvement in complexity also possible for (noiseless) OFFO
algorithms? We answer this question positively in what follows.

The theory developed here combines elements of standard adaptive
regularization methods such as \al{AR$p$}
\cite{BirgGardMartSantToin17} and of the OFFO approaches of
\cite{WuWardBott18} and \cite{GrapStel22}. We exhibit an OFFO
regularization method whose iteration complexity is identical to that
obtained when objective function values are used.  In particular, we
consider convergence to approximate first-order and second-order
critical points, and provide sharp complexity bounds depending on the
degree of derivatives used.

The paper is organized as follows. After introducing the new algorithm
in Section~\ref{algo-s}, we present and discuss a first-order worst-case complexity analysis
in Section~\ref{complexity-s}, while convergence to approximate
second-order minimizers is considered in Section~\ref{2nd-s}. Some
numerical expirements showing the impact of noise are then presented in Section~\ref{numerics_s}.
Conclusions and perspectives are outlined in Section~\ref{concl-s}.

\numsection{An OFFO adaptive regularization algorithm}~\label{algo-s}

\noindent
We now consider the problem of finding approximate minimizers of the
unconstrained nonconvex optimization problem
\beqn{problem}
\min_{x\in \Re^n}f(x),
\eeqn
where $f$ is a sufficiently smooth function from $\Re^n$ into $\Re$.
As motivated in the introduction, our aim is to design an algorithm in
which the objective function value is never computed.
Our approach is based on regularization methods.  In such methods, a model of
the objective function is build by ``regularizing'' a truncated Taylor
expansion of degree $p\geq 1$. We now detail the assumption on the
problems that ensure this approach makes sense.

\noindent
\textbf{AS.1} $f$ is $p$ times continuously differentiable in $\Re^n$. 

\noindent
{\bf AS.2} There exists a constant $f_{\rm low}$ such that
$f(x) \geq f_{\rm low}$ for all $x\in \Re^n$.

\noindent
\textbf{AS.3} The $p$th derivative of $f$ is globally Lipshitz
continuous, that is, there exist a non-negative constant $L_p$ such that
\beqn{LipHessian}
\|\nabla_x^p f(x) - \nabla_x^p f(y)\| \leq L_p \|x-y\| \, \text{ for all } x,y \in \Re^n,
\eeqn
where $\|.\| $ denotes the Euclidean norm in $\Re^n$.

\noindent
\textbf{AS.4}
If $p > 1$, there exists a constant $\khigh\leq 0$ such that
\beqn{Lipitensor}
\min_{\|d\|\leq1} \nabla_x^i f(x) [d]^i \geq \khigh \tim{for all}
x\in\Re^n \tim{and} i\in\iibe{2}{p},
\eeqn
where $\nabla_x^i f(x)$ is the $i$th derivative tensor of
$f$ computed at $x$, and where $T[d]^i$ denotes the
$i$-dimensional tensor $T$ applied on $i$ copies of the vector $d$.
(For notational convenience, we set $\khigh = 0$ if $p=1$.)

We note that AS.4 is weaker than assuming uniform boundedness of the
derivative tensors of degree two and above (there is no upper bound on
the value of $\nabla_x^if(x) [d]^i$), or, equivalently, Lipschitz
continuity of derivatives of degree one to $p-1$.

\subsection{The \al{OFFAR$p$} algorithm}

Adaptive regularization methods are iterative schemes which compute a step from
an iterate $x_k$ to the next by approximately minimizing a $p$th
degree regularized model $m_k(s)$ of $f(x_k+s)$ of the form
\beqn{model}
m_k(s) \eqdef T_{f,p}(x_k,s) + \frac{\sigma_k}{(p+1)!} \|s\|^{p+1},
\eeqn
where $T_{f,p}(x,s)$ is the $p$th order Taylor expansion of
functional $f$ at $x$ truncated at order $p$, that is,
\beqn{taylor model}
T_{f,p}(x,s) \eqdef f(x) + \sum_{i=1}^{p} \frac{1}{i!}\nabla_x^i f(x)[s]^i.
\eeqn
In \req{model}, the $p$th order Taylor series is ``regularized'' by adding the term
$\frac{\sigma_k}{(p+1)!} \|s\|^{p+1}$, where $\sigma_k$ is known as the ``regularization
parameter''. This term guarantees that $m_k(s)$ is bounded below and thus
makes the procedure of finding a step $s_k$ by (approximately) minimizing
$m_k(s)$ well-defined. Our proposed algorithm follows the outline line of existing
\al{AR$p$} regularization methods
\cite{CartGoulToin11d,BirgGardMartSantToin17,CartGoulToin22}, with
the significant difference that the objective function $f(x_k)$ is
never computed, and therefore that the ratio of achieved to predicted
reduction (a standard feature for these methods) cannot be used to
accept or reject a potential new iterate and to update the regularization parameter.  Instead,
such potential iterates are always accepted and the
regularization parameter is updated in a manner independent of this
ratio. We now state the resulting \al{OFFAR$p$} algorithm in detail \vpageref{OFFARp}.

\algo{OFFARp}{ OFFO adaptive regularization of degree $p$ (\tal{OFFAR$p$})}{
	\begin{description}
	   \item[Step 0: Initialization: ] An initial point $x_0\in \Re^n$, a regularization
		parameter $\nu_0>0$ and a requested final gradient accuracy
		$\epsilon_1 \in (0,1]$ are given, as well as the parameters
		\beqn{hyparam}
		 \theta_1 > 1 \tim{ and } \vartheta \in (0,1].
		\eeqn
		Set $k=0$.
	   \item[Step 1: Check for termination: ] Evaluate
                $g_k=\nabla_x^1 f(x_k)$. Terminate with $x_\epsilon = x_k$ if
		\beqn{stopcondhilb}
		\|g_k\| \leq \epsilon_1.
		\eeqn
                Else, evaluate $\{\nabla_x^if(x_k)\}_{i=2}^p$.
	      \item[Step 2: Step calculation: ]
                If $k > 0$, set
                \beqn{muk-def}
                \mu_{1,k} = \frac{p!\|g_k\|}{\|s_{k-1}\|^p}-\theta_1\sigma_{k-1}
                \eeqn
                and select
		\beqn{sigkupdate}
		\sigma_k  \in \left[ \vartheta \nu_k , \max\left( \nu_k,\mu_{1,k}\right) \right].
		\eeqn
                Otherwise (i.e.\ if $k=0$), set $\sigma_0=\mu_0=\nu_0$.\\
                Then compute a step $s_k$  which sufficiently reduces the
                model $m_k$ defined in \req{model} in the sense that 
		\beqn{descent}
		 m_k(s_k) - m_k(0) <  0
		\eeqn
		and
		\beqn{gradstep}
                \|\nabla_s^1 T_{f,p}(x_k,s_k)\| \leq  \theta_1 \frac{\sigma_k}{p!}\|s_k\|^p.
		\eeqn
	   \item[Step 3: Updates. ]
		Set
                \beqn{accept}
                x_{k+1} = x_k + s_k
                \eeqn
                and
		\beqn{vkupdate}
		\nu_{k+1} =  \nu_k + \nu_k \| s_k\|^{p+1}.
		\eeqn
		Increment $k$ by one and go to Step~1.
	\end{description}
}

\noindent
The test \req{gradstep} follows \cite{GratToin21} and extends the more
usual condition where the step $s_k$ is chosen to ensure that
\[
\|\nabla_s^1 m_k(s_k)\| \leq  \theta_1 \|s_k\|^p.
\]
It is indeed easy to verify that \req{gradstep} holds at a local
minimizer of $m_k$ with $\theta_1\geq 1$ (see  \cite{GratToin21} for
details). The motivation for the introduction of $\mu_{1,k}$ in
\req{muk-def} and \req{sigkupdate} will become clearer after Lemma~\ref{useful}.

\numsection{Evaluation complexity for the \tal{OFFAR$p$} algorithm}\label{complexity-s}

\noindent
Before discussing our analysis of evaluation complexity, we first restate some
classical lemmas for \al{AR$p$} algorithms, starting with Lipschitz error bounds.

\llem{lipschitz}{
  Suppose that AS.1 and AS.3 hold. Then
  \beqn{Lip-f}
  |f(x_{k+1})- T_{f,p}(x_k,s_k)| \leq \frac{L_p}{(p+1)!} \|s_k\|^{p+1},
  \eeqn
  and
  \beqn{Lip-g}
  \|g_{k+1}-\nabla_s^1 T_{f,p}(x_k,s_k) \| \leq \frac{L_p}{p!} \|s_k\|^p.
  \eeqn
}

\proof{This is a standard result (see \cite[Lemma~2.1]{CartGoulToin20b}
  for instance).
}

\noindent
We start by stating a simple lower bound on the Taylor series'
decrease.

\llem{model-decrease}{
  \beqn{Tdecr}
  \Delta T_{f,p}(x_k,s_k) \eqdef T_{f,p}(x_k,0)-T_{f,p}(x_k,s_k)
  \geq \frac{\sigma_k}{(p+1)!} \|s_k\|^{p+1}.
  \eeqn
}

\proof{The bound directly results from \req{descent} and \req{model}.}

\noindent
This and AS.2 allow us to establish a lower bound on the decrease in
the objective function (although it is never computed).

\llem{Ifsigmabig}{
Suppose that AS.1 and AS.3 hold and that $\sigma_k \geq 2 L_p$. Then
\beqn{sigma-upper}
f(x_{k}) - f(x_{k+1}) \geq \frac{\sigma_k}{2(p+1)!} \|s_k\|^{p+1}.
\eeqn
}

\proof{
  From \req{Lip-f} and \req{Tdecr}, we obtain that
  \[
  f(x_{k}) - f(x_{k+1}) \geq
  \frac{\sigma_k-L_p}{(p+1)!}\|s_k\|^{p+1}
  \]
  and \req{sigma-upper} immediately follows from our assumption on $\sigma_k$.
} 

\noindent
The next lemma provides a useful lower bound on the step length, in the spirit of
\cite[Lemma~2.3]{BirgGardMartSantToin17} or \cite{GratToin21}.

\llem{useful}{
  Suppose that AS.1 and AS.3 hold. Then, for all $k\geq0$,
  \beqn{crucial}
  \|s_k\|^p > \frac{p!}{L_p + \theta_1\sigma_k } \| g(x_{k+1})\|,
  \eeqn
  \beqn{muk-upper}
  \mu_{1,k} \leq \max[\nu_0, L_p]
  \eeqn
  and
  \beqn{signu}
  \sigma_k \leq L_p + \nu_k.
  \eeqn
}

\proof{
  Successively using the triangle inequality, condition \req{gradstep} and \req{Lip-g}, we deduce that
  \[
  \|{g(x_{k+1})}\| \leq \|g(x_{k+1})- \nabla_s^1 T_{f,p}(x_k,s_k)\| + \|\nabla_s^1 T_{f,p}(x_k,s_k)\|
  \leq \bigfrac{1}{p!} L_p \|s_k\|^p +\theta_1 \frac{\sigma_k}{p!}\|s_k\|^p.
  \]
  The inequality \req{crucial} follows by rearranging the
  terms. Combining this inequality with \req{muk-def} and the identity
  $\mu_0=\nu_0$ then gives \req{muk-upper}, from which \req{signu}
  directly follows using \req{sigkupdate}.
}

\noindent
Observe that \req{muk-upper} motivates
our choice in \req{sigkupdate} to allow the regularization parameter to be of size $\mu_{1,k}$.

Inspired by \cite[Lemma~7]{GrapStel22}, we now establish an upper bound on the number of
iterations needed to enter the algorithm's phase where Lemma~\ref{Ifsigmabig}
applies and thus all iterations produce a decrease in the objective function.

\llem{sigmanotbig}{
Suppose that AS.1 and AS.3 hold, and that the \al{OFFAR$p$} algorithm
does not terminate before or at iteration of index
\beqn{kstardef}
k \geq k_*
\eqdef
\left\lceil\left(\frac{2L_p}{\epsilon_1\vartheta p!}\left((1+\theta_1)\frac{L_p}{\sigma_0}+\theta_1\right)\right)^\sfrac{p+1}{p}\right\rceil.
\eeqn
Then,
\beqn{vkbig}
\nu_k \geq \frac{2L_p}{\vartheta}
\eeqn
which  implies that
\beqn{sigmabig}
\sigma_k \geq 2L_p.
\eeqn
}

\proof{
Note that \eqref{sigmabig} is a direct consequence of \eqref{sigkupdate} if \eqref{vkbig} is true.
Suppose the opposite and that for some $k \geq k_*$, $\nu_k <
\frac{2L_p}{\vartheta}$. Since $\nu_k$ is a non-decreasing sequence, we
have that $\nu_j < \frac{2L_p}{\vartheta}$ for $j \in\iiz{k}$.
Successively using the form of the $\nu_k$ update rule
\eqref{vkupdate}, \req{crucial}, \eqref{signu} and the fact that,
if the algorithm has reached iteration $k_*$, it must
be that \req{stopcondhilb} has failed for all iterations of index at
most $k_*$, we derive that
 \begin{align*}
 \nu_k &> \sum_{j=0}^{k-1} \nu_j \| s_j\|^{p+1} 
 \geq \sum_{j=0}^{k-1} \nu_j \left(\frac{p!\| g(x_{j+1})\|}{L_p+\theta_1\sigma_j}\right)^{\sfrac{p+1}{p}}
 \geq \sum_{j=0}^{k-1} \nu_j \left(\frac{p!\| g(x_{j+1})\|}{L_p+\theta_1 (L_p+\nu_j)}\right)^{\sfrac{p+1}{p}}\\
 &=\sum_{j=0}^{k-1} \nu_j^{-\sfrac{1}{p}} \left(\frac{p!\| g(x_{j+1})\|}{(1+\theta_1)\frac{L_p}{\nu_j}+\theta_1}\right)^{\sfrac{p+1}{p}}
 >\sum_{j=0}^{k-1} \nu_j^{-\sfrac{1}{p}} \left(\frac{p!\| g(x_{j+1})\|}{(1+\theta_1)\frac{L_p}{\sigma_0}+\theta_1}\right)^{\sfrac{p+1}{p}}\\
& >\bigfrac{ k_*\, \vartheta^\sfrac{1}{p} (p! \epsilon_1)^{\sfrac{p+1}{p}}}{(2L_p)^\sfrac{1}{p} \left((1+\theta_1)\frac{L_p}{\sigma_0}+\theta_1\right)^\sfrac{p+1}{p} }.
 \end{align*}
Substituting the definition of $k_*$ in the last inequality, we obtain that 
\[
\frac{2L_p}{\vartheta} < \nu_{k_*} < \frac{2L_p}{\vartheta}, 
\]
which is impossible. Hence no index $k \geq k_*$ exists such that $\nu_k <
\frac{2L_p}{\vartheta}$ and \eqref{vkbig} and \eqref{sigmabig} hold. 
}

\noindent
Observe that \req{kstardef} depends on the ratio
$\frac{L_p}{\sigma_0}$ which is the fraction by which $\sigma_0$
underestimates the Lipschitz constant. This ratio will percolate in
the rest of our analysis. We now define
\beqn{k1def}
k_1 = \min\left\{ k\geq 1 \mid \nu_k \geq \frac{2L_p}{\vartheta}\right\},
\eeqn
the first iterate such that significant objective function decrease is
guaranteed.
The next series of Lemmas provide bounds on $f(x_{k_1})$ and $\sigma_{k_1}$, which in turn will allow
establishing an upper bound on the regularization parameter.
We start by proving an upper bound on $s_k$ generalizing those proposed
in \cite{CartGoulToin11,GratToin21} to the case where $p$ is arbitrary.

\llem{stepgkbound}{
Suppose that AS.1 and AS.4 hold. At each iteration $k$, we have that
\beqn{skbound}
\| s_k\| \leq
2 \eta +2 \left(\frac{(p+1)! \|g_k \|}{\sigma_k}\right)^\sfrac{1}{p},
\eeqn
where 
\beqn{kaphigh}
\eta = \sum_{i=2}^{p}  \left[\frac{\max[0,-\khigh](p+1)!}{i!\,\vartheta \nu_0}\right]^\sfrac{1}{p-i+1}.
\eeqn
}
\proof{If $p=1$, we obtain from \req{descent} and the Cauchy-Schwarz
  inequality that
  \[
  \half \sigma_k \|s_k\|^2 < -g_k^Ts_k \leq \|g_k\|\,\|s_k\|
  \]
  and \req{skbound} holds with $\eta = 0$.
  Suppose now that $p>1$. Again \req{descent} gives that
  \[
  \frac{\sigma_{k}}{(p+1)!} \|s_k\|^{p+1}
  \leq  -g_k^Ts_k -  \sum_{i=2}^p \frac{1}{i!}\nabla_x^i f(x_k)  [ s_k]^i 
  \leq \| g_k\| \| s_k\| +  \sum_{i=2}^p \frac{\max[0,-\khigh]}{i!} \| s_k\|^i. 
  \]
  Applying now the Lagrange bound for polynomial roots \cite[Lecture~VI, Lemma~5]{Yap99}
  with $x=\|s_k\|$, $n=p+1$, $a_0 = 0$, $a_1=\| g_k\|$,
  $a_i = \max[0,-\khigh]/i!$ $i \in \iibe{2}{p}$ and $a_{p+1}
  =\sigma_k/(p+1)!$, we know from \req{descent} that the equation
 $\sum_{i=0}^na_ix^i = 0$ admits at least one strictly positive root,
 and we may thus derive that
 \begin{align*}
 \| s_k\|
  &\leq 2 \left( \frac{ (p+1)!\|g_k\|}{\sigma_k}\right)^\sfrac{1}{p}
   + 2 \sum_{i=2}^{p}  \left[\frac{ \max[0,-\khigh](p+1)!}{i!\,\sigma_k}\right]^\sfrac{1}{p-i+1} \\
  &\leq 2 \left( \frac{ (p+1)!\|g_k\|}{\sigma_k}\right)^\sfrac{1}{p}
   + 2 \sum_{i=2}^{p}  \left[\frac{\max[0,-\khigh](p+1)!}{i!\,\vartheta \nu_k}\right]^\sfrac{1}{p-i+1} \\
  &\leq 2 \left( \frac{ (p+1)!\|g_k\|}{\sigma_k}\right)^\sfrac{1}{p}
   + 2 \sum_{i=2}^{p}  \left[\frac{\max[0,-\khigh](p+1)!}{i!\,\vartheta \nu_0}\right]^\sfrac{1}{p-i+1}, \\
 \end{align*}
and \req{skbound} holds with \req{kaphigh}.
}

\noindent
Our next step is to prove that $\nu_{k_1}$ is bounded by constants only
depending on the problem and the fixed algorithmic parameters.

\llem{vk1bound}{Suppose that AS.1, AS.3 and AS.4 hold. Let $k_1$ be defined by \eqref{k1def}.
We have that, 
\beqn{vk1boundexpr}
\nu_{k_1} \leq \nu_{\max}
= \max \left[
 \sigma_0
 +\sigma_0\left(2 \eta+2\left(\frac{(p+1)!\|g_0\|}{\sigma_0}\right)^\sfrac{1}{p}\right)^{p+1} , \,
 \frac{2\kappa_1L_p}{\vartheta}
      \right]
\eeqn
where $\eta$ is defined in \req{kaphigh} and
\beqn{kap1}
\kappa_1
\eqdef 1+  2^{2p+1}  \eta^{p+1} + 2^{2p+1}
\left[\frac{p+1}{\vartheta}\left((1+\theta_1)\frac{L_p}{\sigma_0}+\theta_1\right)\right]^\sfrac{p+1}{p}.
\eeqn
}
\proof{
If $k_1 = 1$, we have that 
\[
\nu_1 = \sigma_0 + \sigma_0 \|s_0 \|^{p+1}. 
\]
Using Lemma~\ref{stepgkbound} to bound $\| s_0\|^{p+1}$, we derive the
bound corresponding to the first term in the maximum of
\eqref{vk1boundexpr}. Suppose now that $k_1 \geq 2$.
Successively using \eqref{vkupdate}, Lemma~\ref{stepgkbound}, the fact
that $(x+y)^{p+1} \leq 2^{p}(x^{p+1} + y^{p+1})$, the updates rule for
$\nu_k$ \eqref{vkupdate} and $\sigma_k$ \eqref{sigkupdate} and
Lemma~\ref{useful}, we derive that
\begin{align*}
\nu_{k_1} &= \nu_{k_1 - 1} + \nu_{k_1 -1} \|s_{k_1 -1} \|^{p+1} \\
&\leq  \nu_{k_1 - 1} +  \nu_{k_1 -1} \left( 2 \left( (p+1)! \frac{\| g_{k_1-1}\|}{\sigma_{k_1-1}}\right)^\sfrac{1}{p} + 2  \eta\right)^{p+1} \\
&\leq \nu_{k_1 - 1} + 2^{p} \nu_{k_1 - 1} \left[ 2^{p+1}\eta^{p+1}
     + 2^{p+1} \left({\frac{ (p+1)! \|g_{k_1-1} \|}{\sigma_{k_1-1}}}\right)^\sfrac{p+1}{p}  \right] \\
&\leq \nu_{k_1 - 1} + 2^{2p+1} \nu_{k_1 - 1}
      \left[ \eta^{p+1} + \left({\frac{(p+1)! \|g_{k_1-1} \|}{\vartheta \,\nu_{k_1-1}}}\right)^\sfrac{p+1}{p}  \right] \\
&\leq \nu_{k_1 - 1} + 2^{2p+1} \nu_{k_1-1} \eta^{p+1}
     + 2^{2p+1}\left(\frac{(p+1)!}{\vartheta}\right)^\sfrac{p+1}{p}\frac{\|g_{k_1-1} \|^\sfrac{p+1}{p}}{{\nu_{k_1-1}}^\sfrac{1}{p}}\\
&\leq  \nu_{k_1 - 1} +  2^{2p+1} \nu_{k_1-1} \eta^{p+1}
     + 2^{2p+1} \!\left[\frac{(p+1)!}{\vartheta \, p!}(L_p+ \theta_1\sigma_{k_1-2})\right]^\sfrac{p+1}{p}
     \!\!\nu_{k_1-1}^{-\sfrac{1}{p}}\|s_{k_1-2} \|^{p+1}.   \\
\end{align*}
Now $\nu_k$ is a non decreasing sequence, and therefore, using
\req{signu} and the identity $\nu_0=\sigma_0$, 
\begin{align*}
\nu_{k_1}
&\leq \nu_{k_1 - 1} +  2^{2p+1} \nu_{k_1-1} \eta^{p+1} 
       + 2^{2p+1}\! \left[\frac{p+1}{\vartheta}\big((1+\theta_1)L_p+ \theta_1 \nu_{k_1-2}\big)\right]^\sfrac{p+1}{p}
       \!\!\nu_{k_1-2}^{-\sfrac{1}{p}}\|s_{k_1-2} \|^{p+1}  \\
&\leq \nu_{k_1 - 1} +  2^{2p+1} \nu_{k_1-1} \eta^{p+1}
      + 2^{2p+1}\!\left[\frac{p+1}{\vartheta}\left((1+\theta_1)\frac{L_p}{\nu_{k_1-2}}+\theta_1\right)\right]^\sfrac{p+1}{p}
      \!\!\nu_{k_1-2}^{-\sfrac{1}{p}} \nu_{k_1-2}^\sfrac{p+1}{p}\|s_{k_1-2} \|^{p+1} \\
&\leq \nu_{k_1 - 1} +  2^{2p+1} \nu_{k_1-1} \eta^{p+1} 
      + 2^{2p+1} \!\left[\frac{p+1}{\vartheta}\left((1+\theta_1)\frac{L_p}{\sigma_0}+\theta_1\right)\right]^\sfrac{p+1}{p}
      \!\!\nu_{k_1-2} \|s_{k_1-2} \|^{p+1} \\
&\leq \nu_{k_1 - 1} +  2^{2p+1} \nu_{k_1-1} \eta^{p+1}
      + 2^{2p+1}\!\left[\frac{p+1}{\vartheta}\left((1+\theta_1)\frac{L_p}{\sigma_0}+\theta_1\right)\right]^\sfrac{p+1}{p}
      \!\!(\nu_{k_1-1} - \nu_{k_1-2}) \\
&\leq \nu_{k_1 - 1} +  2^{2p+1} \nu_{k_1-1} \eta^{p+1}
      + 2^{2p+1}\!\left[\frac{p+1}{\vartheta}\left((1+\theta_1)\frac{L_p}{\sigma_0}+\theta_1\right)\right]^\sfrac{p+1}{p}
      \!\!\nu_{k_1-1}. \\
\end{align*}
We then obtain the second part of \eqref{vk1boundexpr} by observing
that $ \nu_{k_1-1} \leq \bigfrac{2L_p}{\vartheta}$. 
}

\noindent
This result allows us to establish an upperbound on $f(x_{k_1})$ as a function of $\nu_{\max}$.

\llem{fk1boundsigmak1}{
	Suppose that AS.1, AS.3 and AS.4 hold. Then
	\beqn{fk1bound}
	f(x_{k_1}) \leq f(x_0) + \frac{1}{(p+1)!}\left(\frac{L_p}{\sigma_0} \nu_{\max} +\vartheta \sigma_0\right).
	\eeqn
}
\proof{
	From \eqref{Lip-f} and \eqref{Tdecr}, we know that
	\beqn{gendecr}
	f(x_{j+1}) - f(x_j) \leq (L_p-\sigma_j) \frac{\|s_j\|^{p+1}}{(p+1)!}.
	\eeqn
	Using now the identity $\sigma_0=\nu_0$, \eqref{vkupdate} and
        the fact that $\nu_k$ is a non-decreasing function, we derive
        that 
	\beqn{skbounds}
	\nu_{k_1} \geq \sigma_0 + \sigma_0 \sum_{j=0}^{k_1-1} \|s_j\|^{p+1}.
	\eeqn
	Summing the inequality \eqref{gendecr} for $j \in
        \iiz{k_1-1}$ and using \eqref{skbounds}, \eqref{vkupdate}
        and \eqref{sigkupdate}, we deduce that
	\begin{align*}
	f(x_{k_1}) &\leq f(x_0) + \frac{L_p}{(p+1)!} \sum_{j=0}^{k_1-1}\|s_j\|^{p+1}
          - \frac{1}{(p+1)!} \sum_{j=0}^{k_1-1} \sigma_j \|s_j\|^{p+1} \\
	&\leq  f(x_0) + \frac{L_p}{(p+1)!} \left(\frac{\nu_{k_1} - \sigma_0}{\sigma_0}\right)
          - \frac{1}{(p+1)!} \sum_{j=0}^{k_1-1} \vartheta \nu_j \| s_j\|^{p+1} \\
	&\leq f(x_0) + \frac{L_p}{(p+1)!} \left(\frac{\nu_{k_1} - \sigma_0}{\sigma_0}\right)
          - \frac{\vartheta}{(p+1)!} (\nu_{k_1} - \sigma_0) . \\ 
	\end{align*}
	We then obtain \eqref{fk1bound} by ignoring the negative terms
        in the right-hand side of this last inequality and using
        Lemma~\ref{vk1bound} to bound $\nu_{k_1}$. 
}

\noindent
The two bounds in Lemma ~\ref{fk1boundsigmak1} and
Lemma~\ref{vk1bound} are useful in that they now imply an upper bound on
the regularization parameter, a crucial step in standard theory for
regularization methods.

\llem{uppersigmak}{Suppose that AS.1, AS.3 and AS.4 hold. Suppose also
  that $k \geq k_1$.  Then  
\beqn{boundisgmamax}
\sigma_k
\leq \sigma_{\max}
\eqdef \max\left[\bigfrac{2(p+1)!}{\vartheta}\left(f(x_0)-f_{\rm low}+\frac{1}{(p+1)!}\left(\frac{L_p}{\sigma_0} \nu_{\max}+\vartheta\sigma_0 \right)\right)
+ \nu_{\max}, L_p,\nu_0\right].
\eeqn}
\proof{
Let $j \in \iibe{k_1}{k}$. By the definition of $k_1$ in \eqref{k1def},  $\sigma_j \geq 2L_p$. 
From Lemma~\ref{Ifsigmabig}, we then have that  
\[
f(x_{j}) - f(x_{j+1})
\geq \frac{\sigma_j}{2(p+1)!} \|s_j \|^{p+1}
\geq \vartheta  \frac{\nu_j}{2(p+1)!} \|s_j \|^{p+1}.
\] 
Summing the previous inequality from $j = k_1$ to $k-1$ and using the
$\nu_j$ update rule  \eqref{vkupdate} and AS.2, we deduce that
\[
f(x_{k_1}) - f_{\rm low}\geq f(x_{k_1}) - f(x_k) \geq  \frac{\vartheta}{2(p+1)!} (\nu_k - \nu_{k_1}).
\]
Rearranging the previous inequality and using Lemma~\ref{vk1bound},
\beqn{nuk-upper}
\nu_k \leq \frac{2(p+1)!}{\vartheta} \left( f(x_{k_1}) - f_{\rm low} \right) + \nu_{\max}.
\eeqn
Combining now Lemma~\ref{fk1boundsigmak1} (to bound $f(x_{k_1})$),
\req{sigkupdate} and \req{muk-upper} gives \eqref{boundisgmamax}.
}

\noindent
We may now resort to the standard ``telescoping sum'' argument to obtain the
desired evaluation complexity bound.

\lthm{complexity}{
  Suppose that AS.1--AS.4 hold. Then the \al{OFFAR$p$} algorithm requires at most
  \[
  \left[\kap{OFFARp}\,\left( f(x_0)-f_{\rm low}
    + \frac{1}{(p+1)!}\left(\frac{L_p}{\sigma_0}\nu_{\max}+\vartheta\sigma_0\right)\right)
   + \left(\frac{2L_p}{\vartheta p!}\left(\frac{L_p}{\sigma_0}+\theta_1\right)\right)^\sfrac{p+1}{p}\right]
\epsilon_1^{-\sfrac{p+1}{p}}+2
  \]
  iterations and evaluations of $\{\nabla_x^i f\}_{i = 1}^p$  to
  produce a vector $x_\epsilon\in \Re^n$ such that
  $\| {g(x_\epsilon)} \| \leq \epsilon_1$, where
  \[
  \kap{OFFARp}
  \eqdef 2(p+1)!\,\sigma_{\max}^{1/p}\left(\frac{1}
         {\vartheta p!}\left(\frac{L_p}{\sigma_0} + \vartheta\theta_1\right)\right)^\sfrac{p+1}{p}
  \]
where $\sigma_{\max}$ is defined in Lemma~\ref{uppersigmak} and
$\nu_{\max}$ is defined in Lemma~\ref{vk1bound}.
}

\proof{
Suppose that the algorithm terminates at an iteration $k < k_1$, where
$k_1$ is given by \req{k1def}. The desired conclusion then follows
from the fact that, by this definition and Lemma~\ref{sigmanotbig},
\beqn{k1bound}
k_1\leq k_* \leq
\left(\frac{2L_p}{\epsilon_1 \vartheta p!}\left(\frac{L_p}{\sigma_0}+\theta_1\right)\right)^\sfrac{p+1}{p} +1.
\eeqn
Suppose now that the algorithm has not terminated at iteration $k_1$
and consider an iteration $j \geq k_1$. From $k_1$ definition \eqref{k1def} and Lemma~\ref{uppersigmak}, we have that $2L_p\leq \sigma_j \leq \sigma_{\max} $. 
Since $\sigma_j \geq 2L_p$, Lemma~\ref{Ifsigmabig} is valid for iteration $j$.
But $\sigma_j\in[\vartheta \sigma_0,\sigma_{\max}]$ because of Lemma~\ref{uppersigmak} and
$\|g(x_{j+1})\|\geq\epsilon_1$ before termination, and we therefore deduce that
\beqn{inproof1}
f(x_j) - f(x_{j+1})
\geq \frac{\sigma_j \|s_j \|^{p+1}}{2(p+1)!}
\geq \frac{\sigma_j (p!)^\sfrac{p+1}{p} \|g(x_{j+1})\|^\sfrac{p+1}{p}}
          { 2 (p+1)!(L_p + \theta_1\sigma_j)^\sfrac{p+1}{p}}  
\geq\frac{(p!)^\sfrac{p+1}{p} \epsilon_1^\sfrac{p+1}{p}}
          {2(p+1)!\sigma_{\max}^\sfrac{1}{p}\left(\frac{L_p}{\vartheta \sigma_0} + \theta_1\right)^\sfrac{p+1}{p}}.  
\eeqn
Summing this inequality from $k_1$ to $k\geq k_1$ and using AS.3, we
obtain that
\beqn{inproof2}
f(x_{k_1}) - f_{\rm low}\geq f(x_{k_1}) - f(x_k) \geq  \frac{(k-k_1)}{\kap{OFFARp}} \epsilon_1^\sfrac{p+1}{p}.
\eeqn
Rearranging the terms of the last inequality and using \req{k1bound}
and Lemma~\ref{fk1boundsigmak1} then yields the desired result.
} 

\noindent
While this theorem covers all model's degrees, it is worthwhile to
isolate the most commonly used cases.

\lcor{complexity2}{
  Suppose that AS.1--AS.3 hold and that $p=1$. Then the \al{OFFAR1} algorithm requires at most
  \[
  \left[\frac{4\sigma_{\max}}{\vartheta}\left(\frac{L_1}{\sigma_0} + \vartheta \theta_1\right)^2
    \!\left[ f(x_0)\!-\!f_{\rm low}
    \!+ \!\frac{1}{2}\left(\frac{L_1}{\sigma_0} \nu_{\max}+\vartheta\sigma_0\right)\right]
   \! + \!\left(\frac{2L_1}{\vartheta}\left(\frac{L_1}{\sigma_0}+\theta_1\right)\right)^2\right]
  \epsilon_1^{-2}+2
  \]
  iterations and evaluations of the gradient to produce a vector $x_\epsilon\in \Re^n$ such that
  $\| {g(x_\epsilon)} \| \leq \epsilon_1$, 
  where $\sigma_{\max}$ is defined in Lemma~\ref{uppersigmak} and
  $\nu_{\max}$ is defined in Lemma~\ref{vk1bound}. If $p=2$ and AS.4
  holds, the \al{OFFAR2} algorithm requires at most 
  \[
  \left[\frac{12\sigma_{\max}^{1/2}}{(2\vartheta)^\sfrac{3}{2}}
   \! \left(\frac{L_2}{\sigma_0}+ \vartheta\theta_1\right)^\sfrac{3}{2}\!\left[ f(x_0)\!-\!f_{\rm low}
   \! + \! \frac{1}{6}\left(\frac{L_2}{\sigma_0} \nu_{\max}+\vartheta\sigma_0\right)\right]
   \! + \! \left(\frac{L_2}{\vartheta}\left(\frac{L_2}{\sigma_0}+\theta_1\right)\right)^\sfrac{3}{2}\right]
  \epsilon_1^{-\sfrac{3}{2}}+2
  \]
  iterations and evaluations of the gradient and Hessian to achieve
  the same result.
}

\noindent
We now prove that the complexity bound stated by
Theorem~\ref{complexity} is sharp in order.

\lthm{sharpness}{
Let $\epsilon_1\in(0,1]$ and $p \geq 1$.  Then there exists a $p$ times continuously
differentiable function $f_p$ from $\Re$ into $\Re$ such that the
\al{OFFAR$p$} applied to $f_p$ starting from the origin takes exactly
$k_\epsilon = \lceil\epsilon_1^{-\sfrac{p+1}{p}}\rceil$ iterations and derivative's
evaluations to produce an iterate $x_{k_\epsilon}$ such that $|\nabla_x^1f_p(x_{k_\epsilon})|\leq\epsilon_1$.
}

\proof{
To prove this result, we first define a sequence of function
and derivatives' values such that the gradients converge sufficiently
slowly and then show that these sequences can be generated by the
\al{OFFAR$p$} algorithm and also that there exists a function $f_p$
satisfying AS.1--AS.4 which interpolate them.

First select $\vartheta = 1$ (implying that $\sigma_k=\nu_k$ for all
$k$), some $\sigma_0=\nu_0>0$ and define, for all $k\in\iiz{k_\epsilon}$,
\beqn{omegak-ex}
\omega_k = \epsilon_1\,\frac{k_\epsilon-k}{k_\epsilon} \in [0,\epsilon_1]
\eeqn
and
\beqn{gk-ex}
g_k = -(\epsilon_1+\omega_k)
\tim{ and }
D_{i,k} = 0,  \ms (i=2,\ldots,p),
\eeqn
so that
\beqn{grange-ex}
|g_k| \in[ \epsilon_1,2\epsilon_1] \subset [0,2] \tim{for all} k\in\iiz{k_\epsilon}.
\eeqn
We then set, for all $k\in\iiz{k_\epsilon}$,
\beqn{sk-ex}
s_k = \left(\frac{p!|g_k|}{\sigma_k}\right)^\frac{1}{p},
\eeqn
so that
\begin{align}
\sigma_k
&\eqdef \sigma_0+\sum_{j=0}^{k-1}\sigma_j|s_j|^{p+1}\label{sig-ex}\\
&= \sigma_0+\sum_{j=0}^{k-1}\sigma_j\left(\frac{p!|g_j|}{\sigma_j}\right)^\frac{p+1}{p}
= \sigma_0+(p!)^\frac{p+1}{p}\sum_{j=0}^{k-1}\frac{(\epsilon_1+\omega_j)^\frac{p+1}{p}}{\sigma_j^\frac{1}{p}}\nonumber\\
&\leq \sigma_0+\left(\frac{(2p!)^{p+1}}{\sigma_0}\right)^\frac{1}{p}\sum_{j=0}^{k-1}\epsilon_1^\frac{p+1}{p}
\leq \sigma_0+\left(\frac{(2p!)^{p+1}}{\sigma_0}\right)^\frac{1}{p}k_\epsilon\epsilon_1^\frac{p+1}{p}
\leq \sigma_0+2\left(\frac{(2p!)^{p+1}}{\sigma_0}\right)^\frac{1}{p}\nonumber\\
& \eqdef \sigma_{\max},\nonumber
\end{align}
where we successively used \req{sk-ex}, \req{gk-ex}, \req{omegak-ex} and
the definition of $k_\epsilon$.
We finally set
\[
f_0 =2^\frac{2p+1}{p}\left(\frac{p!}{\sigma_0}\right)^\frac{1}{p}
\tim{and}
f_{k+1} \eqdef f_k + g_ks_k + \sum_{i=2}^p\frac{1}{i!}D_{i,k}[s_k]^i
= f_k - \left(\frac{p!}{\sigma_k}\right)^\frac{1}{p}(\epsilon_1+\omega_k)^\frac{p+1}{p},
\]
yielding, using \req{sig-ex} and the definition of $k_\epsilon$, that
\[
f_0-f_{k_\epsilon}  =
\sum_{k=0}^{k_\epsilon-1}\left(\frac{p!}{\sigma_k}\right)^\frac{1}{p}(\epsilon_1+\omega_k)^\frac{p+1}{p}
\leq 2^\frac{p+1}{p}\left(\frac{p!}{\sigma_0}\right)^\frac{1}{p} k_\epsilon\epsilon_1^\frac{p+1}{p}
\leq 2^\frac{2p+1}{p}\left(\frac{p!}{\sigma_0}\right)^\frac{1}{p}
= f_0.
\]
As a consequence
\beqn{frange-ex}
f_k\in[0,f_0] \tim{for all} k\in\iiz{k_\epsilon}.
\eeqn
Observe that \req{sk-ex} satisfies \req{descent} (for the model
\req{model}) and \req{gradstep} for $\theta_1 =1$.  Moreover \req{sig-ex} is
the same as \req{vkupdate}-\req{sigkupdate}. Hence the sequence
$\{x_k\}$ generated by 
\[
x_0 = 0
\tim{and}
x_{k+1}=x_k+s_k
\]
may be viewed as produced by the \al{OFFAR$p$} algorithm given
\req{gk-ex}. Defining
\[
T_{k,p}(s) = f_k + g_ks + \sum_{i=2}^p \frac{1}{i!}D_{i,k}s^i,
\]
observe also that
\beqn{fLip-ex}
|f_{k+1}-T_{k,p}(s_k)| = |g_ks_k|
\leq 2(p!)^\frac{1}{p}\sigma_{\max}\left(\frac{\epsilon_1+\omega_k}{\sigma_k}\right)^\frac{p+1}{p}
\leq \frac{2\sigma_{\max}}{p!}\,|s_k|^{p+1}
\eeqn
and
\beqn{gLip-ex}
|g_{k+1}-\nabla_s^1T_{k,p}(s_k)| = |g_{k+1}-g_k|
\leq |\omega_k-\omega_{k+1}|
= \frac{\epsilon_1}{k_\epsilon}
\leq \epsilon_1^\frac{2p+1}{p}
\leq \frac{\sigma_{\max}}{\sigma_k}(\epsilon_1+\omega_k)
= \frac{\sigma_{\max}}{p!}|s_k|^p
\eeqn
(we used $k_\epsilon \leq \epsilon_1^{-\sfrac{p+1}{p}}+1$ and
$\epsilon_1\leq 1$), while, if $p>1$,
\beqn{HLip-ex}
|D_{i,k+1}-\nabla_s^iT_{k,p}(s_k)| = |D_{i,k+1}-D_{i,k}| = 0 \leq \frac{\sigma_{\max}}{p!}\,|s_k|^{p+1-i}
\eeqn
for $i=2,\ldots,p$.
In view of \req{grange-ex}, \req{frange-ex} and
\req{fLip-ex}-\req{HLip-ex}, we may then apply classical Hermite
interpolation to the data given by
$\{(x_k,f_k,g_k,D_{2,k},\ldots,D_{p,k})\}_{k=0}^{k_\epsilon}$ (see
\cite[Theorem~A.9.2]{CartGoulToin22} with
$
\kappa_f = \max[ 2, f_0, 2\sigma_{\max}/p!  ],
$
for instance) and deduce that there
exists a $p$ times continuously differentiable piecewise polynomial function $f_p$
satisfying AS.1--AS.4 and such that, for $k\in\iiz{k_\epsilon}$,
\[
f_k = f_p(x_k),
\ms
g_k = \nabla_x^1 f_p(x_k)
\tim{and}
D_{i,k} = \nabla_x^i f_p(x_k), \ms (i=2,\ldots,p).
\]
The sequence $\{x_k\}$ may thus be interpreted as being produced by
the \al{OFFAR$p$} algorithm applied to $f_p$ starting from $x_0=0$.
The desired conclusion then follows by observing that, from
\req{omegak-ex} and \req{gk-ex},
\[
|g_k| > \epsilon_1 \tim{for} k\in\iiz{k_\epsilon-1}
\tim{and}
|g_{k_\epsilon}| = \epsilon_1.
\]
} 

\noindent
It is remarkable that the complexity bound stated by
Theorems~\ref{complexity} and \ref{complexity-2nd} are identical (in order) to that known for
the standard setting where the objective function is evaluated at each
iteration.  Moreover, the $\calO(\epsilon^{-3/2})$ bound for $p=2$ was
shown in \cite{CartGoulToin18a} to be optimal within a very large class
of second-order methods. One then concludes that, from the sole
viewpoint of evaluation complexity, the computation of the objective
function's values is an unnecessary effort for achieving convergence
at optimal speed.

One may also ask, at this point, if keeping track of $\nu_k$ is
necessary, that is, when considering \al{OFFAR2}, if a simplified
update of the form
\beqn{sigbis}
\sigma_k = \max[\sigma_{k-1},\mu_{1,k}]
\eeqn
would not be sufficient to ensure convergence at the desired rate.  As
we show in appendix, this is not the case, because $\mu_{1,k}$ only
measures change in second derivatives along the direction $s_{k-1}$,
thereby producing an underestimate of $L_p$. As a result, $\sigma_k$
may fail to reach this value and a simplified \al{OFFAR2} algorithm using
\req{sigbis} instead of \req{sigkupdate} may fail to converge altogether.
Another mechanism (such as that provided by $\nu_k$) is thus
necessary to force the growth of the regularization
parameter beyond the Lispchitz constant.

\numsection{Second-order optimality}\label{2nd-s}

If second-derivatives are available and $p \geq 2$, it is also
possible to modify the \al{OFFAR$p$} algorithm to obtain second-order
optimality guarantees.  We thus assume in this section that $p\geq2$
and restate  the algorithm \vpageref{MOFFARp}.

\algo{MOFFARp}{Modified OFFO adaptive regularization of degree $p$ (\al{MOFFAR$p$})}{
	\begin{description}
	   \item[Step 0: Initialization: ] An initial point $x_0\in \Re^n$, a regularization
		parameter $\nu_0>0$, a requested final gradient accuracy
		$\epsilon_1 \in (0,1]$ and a requested final curvature accuracy
		$\epsilon_2 \in (0,1]$ are given, as well as the parameters
		\beqn{hyparam-m}
		\theta_1,\theta_2 > 1 \tim{ and } \vartheta \in (0,1]
		\eeqn
		Set $k=0$.
           \item[Step 1: Check for termination: ] Evaluate $g_k=\nabla_x^1 f(x_k)$ and $\nabla_x^2f(x_k)$. 
                Terminate with $x_\epsilon = x_k$ if
		\beqn{stopcondhilb-m}
		\|g_k\| \leq \epsilon_1
                \tim{ and }
                \lambda_{\min}[\nabla_x^2f(x_k)] \geq -\epsilon_2.
		\eeqn
                Else, evaluate $\{\nabla_x^if(x_k)\}_{i=3}^p$.
	   \item[Step 2: Step calculation: ] If $k > 0$, set
                \beqn{muk-def-m}
                \mu_{1,k} =
                \frac{p!\|g_k\|}{\|s_{k-1}\|^p}-\theta_1\sigma_{k-1},
                \ms
                \mu_{2,k} =
                \frac{(p-1)!\max\Big[0,-\lambda_{\min}[\nabla_x^2 f(x_k)]\Big]}{\|s_{k-1}\|^{p-1}}
                 - \theta_2 \sigma_{k-1}
                \eeqn
                and select
		\beqn{sigkupdate-m}
		\sigma_k  \in \left[ \vartheta \nu_k , \max\left( \nu_k,\mu_{1,k},\mu_{2,k}\right) \right].
		\eeqn
                Otherwise (i.e.\ if $k=0$), set $\sigma_0=\mu_{1,0}=\mu_{2,0}=\nu_0$.\\
                Then compute a step $s_k$  which
		sufficiently reduces the model $m_k$ defined in \req{model} in the sense that
		\beqn{descent-m}
		 m_k(s_k) - m_k(0) <  0,
		\eeqn
		\beqn{gradstep-m}
                \|\nabla_s^1 T_{f,p}(x_k,s_k)\| \leq \theta_1\frac{\sigma_k}{p!}\|s_k\|^p
                \eeqn
                and
                \beqn{curvstep-m}
                \lambda_{\min}[\nabla_s^2 T_{f,p}(x_k,s_k)] \geq  -\theta_2\frac{\sigma_k}{(p-1)!}\|s_k\|^{p-1}.
		\eeqn
	   \item[Step 3: Updates. ]
		Set
                \beqn{accept-m}
                x_{k+1} = x_k + s_k,
                \eeqn
                and
		\beqn{vkupdate-m}
		\nu_{k+1} =  \nu_k + \nu_k \| s_k\|^{p+1}.
		\eeqn
		Increment $k$ by one and go to Step~1.
	\end{description}
}

\noindent
The modified algorithm only differs from that of page~\pageref{OFFARp}
by the addition of the second part of \req{stopcondhilb-m}, the
inclusion of $\mu_{2,k}$ (whose purpose parallels that of $\mu_{1,k}$)
and condition \req{curvstep-m} on the step $s_k$.  As was
the case for \req{gradstep}/\req{gradstep-m}, note that
\req{curvstep-m} holds with $\theta_2=1$ at a second-order
minimizer of the model $m_k(s)$, and is thus achievable for
$\theta_2>1$. Moreover, because the modified algorithm subsumes the original
one, all properties derived in the previous section continue to hold.
In addition, we may complete the
bounds of Lemma~\ref{lipschitz} by noting that AS.3 for $p>1$ also implies that
\beqn{Lip-H}
\|\nabla_x^2 f(x_{k+1})-\nabla_s^2 T_{f,p}(x_k,s_k) \| \leq \frac{L_p}{(p-1)!} \|s_k\|^{p-1}.
\eeqn
We now derive a second-order analog of the step lower bound of Lemma~\ref{useful}.

\llem{useful2}{
  Suppose that AS.1 and AS.3 hold and that the modified algorithm is
  applied. Then, for all $k\geq 0$,
  \beqn{crucial2}
  \|s_k\|^{p-1} > \frac{(p-1)!}{L_p + \theta_2\sigma_k } \max\Big[0,-\lambda_{\min}[\nabla_x^2 f(x_{k+1})]\Big]
  \eeqn
  and
  \beqn{muk-upper-m}
  \mu_{2,k} \leq \max[\nu_0,L_p].
  \eeqn
}

\proof{
  Successively using the triangle inequality, \req{Lip-H} and
  \req{curvstep-m}, we obtain that
  \begin{align*}
  \lambda_{\min}[\nabla_x^2 f(x_{k+1})]
  &= \min_{\|d\|\leq 1} \nabla_x^2 f(x_{k+1})[d]^2\\
  &= \min_{\|d\|\leq 1} \Big[\nabla_x^2 f(x_{k+1})[d]^2 - \nabla_s^2 T_{f,p}(x_k,s_k)[d]^2 + \nabla_s^2 T_{f,p}(x_k,s_k)[d]^2\Big]\\
  &\geq \min_{\|d\|\leq 1} \Big[ \nabla_x^2 f(x_{k+1})[d]^2 - \nabla_s^2 T_{f,p}(x_k,s_k)[d]^2\Big] + \min_{\|d\|\leq 1} \nabla_s^2 T_{f,p}(x_k,s_k)[d]^2\\
  &= \min_{\|d\|\leq 1} \Big[ (\nabla_x^2 f(x_{k+1}) - \nabla_s^2 T_{f,p}(x_k,s_k) )[d]^2\Big] + \lambda_{\min}[\nabla_s^2 T_{f,p}(x_k,s_k)]\\
  &\geq - \|\nabla_x^2 f(x_{k+1}) - \nabla_s^2 T_{f,p}(x_k,s_k)\| - \theta_2\frac{\sigma_k}{(p-1)!}\|s_k\|^{p-1}\\
  &= - \frac{L_p}{(p-1)!}\|s_k\|^{p-1} - \theta_2\frac{\sigma_k}{(p-1)!}\|s_k\|^{p-1},
  \end{align*}
  which proves \req{crucial2}. The bound \req{muk-upper-m} then results
  from the identity $\mu_{2,0} = \nu_0$, \req{sigkupdate-m} and \req{muk-upper-m}.
}

\noindent
Observe that \req{sigkupdate-m}, \req{muk-upper} and \req{muk-upper-m}
ensure that \req{signu} continues to hold.

We now have to adapt our argument since the termination test
\req{stopcondhilb-m} may fail if either its first or its second part
fails.  Lemma~\ref{useful} then gives a lower bound on the step
if the first part fails, while we have to use Lemma~\ref{useful2} otherwise.
This is formalized in the following lemma.

\llem{sigmanotbig-m}{
Suppose that AS.1 and AS.3 hold, and that the \al{MOFFAR$p$} algorithm
has reached iteration of index
\beqn{ks2def}
k \geq k_{**}
\eqdef
\left\lceil
\frac{2L_p}{\kap{both}^{p+1}\vartheta}
 \max\left[\left(\frac{2L_p}{\vartheta}\right)^{\sfrac{1}{p}},\left(\frac{2L_p}{\vartheta}\right)^{\sfrac{2}{p-1}}\right]
 \max\Big[\epsilon_1^{-\sfrac{p+1}{p}},\epsilon_2^{-\sfrac{p+1}{p-1}}\Big]
\right\rceil,
\eeqn
where
\beqn{kapbdef}
\kap{both} \eqdef\min\left[
  \left(\frac{p!}{(1+\theta_1)\frac{L_p}{\sigma_0} + \theta_1 }\right)^\sfrac{1}{p},
  \left(\frac{(p-1)!}{(1+\theta_2)\frac{L_p}{\sigma_0} + \theta_2 } \right)^\sfrac{1}{p-1}
  \right].
\eeqn
Then
\beqn{vkbig-m}
\nu_k \geq \frac{2L_p}{\vartheta},
\eeqn
which  implies that
\beqn{sigmabig-m}
 \sigma_k \geq 2L_p.
\eeqn
}

\proof{
As in Lemma~\ref{sigmanotbig}, \eqref{sigmabig-m} is a direct
consequence of \eqref{sigkupdate-m} if \eqref{vkbig-m} is true.
In order to adapt the
proof of Lemma~\ref{sigmanotbig}, we observe that, at iteration $k$,
\req{crucial} and \req{crucial2} hold and
\[
\|s_k\| > \min\left[
  \left(\frac{p!}{L_p + \theta_1\sigma_k } \| g(x_{k+1})\|\right)^\sfrac{1}{p},
  \left(\frac{(p-1)!}{L_p + \theta_2\sigma_k } \max\Big[0,-\lambda_{\min}[\nabla_x^2 f(x_{k+1})]\Big]\right)^\sfrac{1}{p-1}
  \right]
\]
which, given \req{signu}, that termination has not yet occured and
that $\nu_k > \nu_0=\sigma_0$, implies that
\begin{align}
\|s_k\|&  > \min\left[
  \nu_k^{-\sfrac{1}{p}}\left(\frac{p!}{(1+\theta_1)\frac{L_p}{\sigma_0} + \theta_1 }\right)^\sfrac{1}{p},
  \nu_k^{-\sfrac{1}{p-1}}\left(\frac{(p-1)!}{(1+\theta_2)\frac{L_p}{\sigma_0} + \theta_2 } \right)^\sfrac{1}{p-1}
  \right]
\min\Big[\epsilon_1^\sfrac{1}{p},\epsilon_2^\sfrac{1}{p-1}\Big]\nonumber\\
& \geq \kap{both}\min\left[\nu_k^{-\sfrac{1}{p}},\nu_k^{-\sfrac{1}{p-1}}\right]
\min\Big[\epsilon_1^\sfrac{1}{p},\epsilon_2^\sfrac{1}{p-1}\Big].\label{useful12}
\end{align}
Suppose now that \req{vkbig-m} fails, i.e.\ that for some $k \geq k_{**}$, $\nu_k <
\frac{2L_p}{\vartheta}$. Since $\nu_k$ is a non-decreasing sequence, we
have that $\nu_j < \frac{2L_p}{\vartheta}$ for $j \in\iiz{k}$.
Successively using \eqref{vkupdate-m} and \req{useful12}, we
obtain that
 \begin{align*}
 \nu_k &> \sum_{j=0}^{k-1} \nu_j \| s_j\|^{p+1} 
 \geq \sum_{j=0}^{k-1} \kap{both}^{p+1}
 \min\left[\nu_j^{-\sfrac{1}{p}},\nu_j^{-\sfrac{2}{p-1}}\right]\min\Big[\epsilon_1^\sfrac{1}{p},\epsilon_2^\sfrac{1}{p-1}\Big]^{p+1}\\
 &\geq \sum_{j=0}^{k-1} \kap{both}^{p+1}
 \min\left[\left(\frac{2L_p}{\vartheta}\right)^{-\sfrac{1}{p}},\left(\frac{2L_p}{\vartheta}\right)^{-\sfrac{2}{p-1}}\right]
 \min\Big[\epsilon_1^\sfrac{1}{p},\epsilon_2^\sfrac{1}{p-1}\Big]^{p+1}\\
& = k_{**} \kap{both}^{p+1}
 \min\left[\left(\frac{2L_p}{\vartheta}\right)^{-\sfrac{1}{p}},\left(\frac{2L_p}{\vartheta}\right)^{-\sfrac{2}{p-1}}\right]
 \min\Big[\epsilon_1^\sfrac{1}{p},\epsilon_2^\sfrac{1}{p-1}\Big]^{p+1}.
 \end{align*}
Using the definition of $k_{**}$ in the last inequality, we see that 
\[
\frac{2L_p}{\vartheta} < \nu_{k_{**}} < \frac{2L_p}{\vartheta}, 
\]
which is impossible. Hence no index $k \geq k_{**}$ exists such that $\nu_k <
\frac{2L_p}{\vartheta}$ and \eqref{vkbig-m} and \eqref{sigmabig-m} hold. 
}

\noindent
We then continue to use the theory of the previous section with a value
of $k_1$ now satisfying the improved bound
\beqn{k*k**}
k_1 \leq k_{**},
\eeqn
instead of $k_1\leq k_*$. This directly leads us to the following strengthened complexity result.

\lthm{complexity-2nd}{
  Suppose that AS.1--AS.4 hold and that $p>1$ Then the \al{MOFFAR$p$} algorithm requires at most
  \begin{align*}
    &\left[
      \kap{MOFFARp}\,\left( f(x_0)-f_{\rm low}
    + \frac{1}{(p+1)!}\left(\frac{L_p}{\sigma_0}\nu_{\max}+\vartheta\sigma_0\right)\right)\right]\\
    &\hspace*{20mm} + \frac{2L_p}{\kap{both}^{p+1}\vartheta} 
  \max\left[\left(\frac{2L_p}{\vartheta}\right)^{\sfrac{1}{p}},\left(\frac{2L_p}{\vartheta}\right)^{\sfrac{2}{p-1}}\right]
\max\Big[\epsilon_1^{-\sfrac{p+1}{p}},\epsilon_2^{-\sfrac{p+1}{p-1}}\Big]
     +2
  \end{align*}
  iterations and evaluations of $\{\nabla_x^i f\}_{i = 1}^p$  to produce a vector $x_\epsilon\in \Re^n$ such that
  $\| {g(x_\epsilon)} \| \leq \epsilon_1$ and $\lambda_{\min}[\nabla_x^2f(x_\epsilon)] \geq -\epsilon_2$, where
  \[
  \kap{MOFFARp}
  \eqdef 2(p+1)!\max\left[
     \sigma_{\max}^{1/p}\left(\frac{\frac{L_p}{\sigma_0} + \vartheta \theta_1}{\vartheta p!}\right)^\sfrac{p+1}{p},
     \sigma_{\max}^{2/p-1}\left(\frac{\frac{L_p}{\sigma_0} + \vartheta \theta_2}{\vartheta (p-1)!}\right)^\sfrac{p+1}{p-1}
          \right]
  \]
and where $\sigma_{\max}$ is defined in Lemma~\ref{uppersigmak},
$\nu_{\max}$ is defined in Lemma~\ref{vk1bound} and $\kap{both}$ in \req{kapbdef}.
}

\proof{
  The bound of Theorem~\ref{complexity} remains valid for obtaining a vector $x_\epsilon\in \Re^n$ such that
  $\| {g(x_\epsilon)} \| \leq \epsilon_1$, but we are now interested
  to satisfy the second part of \req{stopcondhilb-m} as well.  Using 
  \req{crucial2} instead of \req{crucial}, we deduce (in parallel to
  \req{inproof1}) that before termination, 
\begin{align*}
f(x_j) - f(x_{j+1})
&\geq \frac{\sigma_j \|s_j \|^{p+1}}{2(p+1)!}\\
&\geq \frac{\sigma_j ((p-1)!)^{\sfrac{p+1}{p-1}}\max[0,-\lambda_{\min}[\nabla_x^2 f(x_{k+1})]]^{\sfrac{p+1}{p-1}}}
           {2(p+1)!(L_p + \theta_2\sigma_j)^\sfrac{p+1}{p-1}}  \\
&\geq \frac{((p-1)!)^\sfrac{p+1}{p-1} \epsilon_2^\sfrac{p+1}{p-1}}
           {2(p+1)!\sigma_{\max}^\sfrac{2}{p-1}\left(\frac{L_p}{\vartheta \sigma_0} + \theta_2\right)^\sfrac{p+1}{p}},  
\end{align*}
so that, summing this inequality from $k_1$ to $k\geq k_1$ and using
AS.3 now gives (in parallel to \req{inproof2}) that, before the second part of
\req{stopcondhilb-m} is satisfied,
\[
f(x_{k_1}) - f_{\rm low}
\geq f(x_{k_1}) - f(x_k)
\geq  \frac{(k-k_1)}{\kap{2nd}} \epsilon_2^\sfrac{p+1}{p-1}
\]
where
\[
\kap{2nd}
\eqdef 2(p+1)!\,\sigma_{\max}^{2/p-1}
\left(\frac{\frac{L_p}{\sigma_0} + \vartheta \theta_2}{\vartheta (p-1)!}\right)^\sfrac{p+1}{p-1}.
           \]
As a consequence, we deduce, using \req{k*k**}, that the second part of
\req{stopcondhilb-m} must hold at the latest after
  \[
  \left[\kap{2nd}\,\left( f(x_0)-f_{\rm low}
    + \frac{1}{(p+1)!}\left(\frac{L_p}{\sigma_0}\nu_{\max}+\vartheta\sigma_0\right)\right)\right] \epsilon_2^{-\sfrac{p+1}{p-1}}
    \!\!+ k_{**} +2
  \]
  iterations and evaluations of the derivatives, where $k_{**}$ is defined
  in \req{ks2def}. Combining this result
  with that of Theorem~\ref{complexity} then yields the desired conclusion.
}

Focusing again on the case where $p=2$ and upperbounding complicated
constants, we may state the following corollary.

\lcor{complexity2-2nd}{
  Suppose that AS.1--AS.4 hold and that $p=2$. Then there exists
  constants $\kappa_*$ such that the \al{MOFFAR2} algorithm requires at most
  \[
  \kappa_* \max\left[\epsilon_1^{-3/2},\epsilon_2^{-3}\right] 
  \]
  iterations and evaluations of the gradient and Hessian to produce a vector $x_\epsilon\in \Re^n$ such that
  $\| {g(x_\epsilon)} \| \leq \epsilon_1$ and $\lambda_{\min}[\nabla_x^2f(x_{k_\epsilon})] \geq -\epsilon_2$.
}

\noindent
We finally prove that the complexity for reaching approximate second
order points, as stated by Theorem~\ref{complexity-2nd}, is also sharp.

\lthm{sharpness2}{
Let $\epsilon_1,\epsilon_2\in(0,1]$ and $p>1$.  Then there exists a $p$ times continuously
differentiable function $f_p$ from $\Re$ into $\Re$ such that the 
\al{MOFFAR$p$} applied to $f_p$ starting from the origin takes exactly
$k_\epsilon = \lceil\epsilon_2^{-\sfrac{p+1}{p-1}}\rceil$ iterations and derivative's
evaluations to produce an iterate $x_{k_\epsilon}$ such that
$|\nabla_x^1f_p(x_{k_\epsilon})|\leq \epsilon_1$
and $\lambda_{\min}[\nabla_x^2f(x_{k_\epsilon})] \geq -\epsilon_2$.
}

\proof{The proof is very similar to that of Theorem~\ref{sharpness}, this
  time taking a uniformly zero gradient but a minimal eigenvalue of
  the Hessian slowly converging to $-\epsilon_2$ from below.  It is
  detailed in appendix.
}

\numsection{The effect of noise}\label{numerics_s}

We have mentioned in the introduction that OFFO algorithms like that
presented above are interesting not only because of their remarkable
theoretical properties covered in the previous sections, but also
because they show remarkable insensitivity to noise\footnote{A similar
behaviour was observed in \cite{GratJeraToin22b} for first-order OFFO methods.}. This section is
devoted to illustrating this statement while, at the same time,
proposing a more detailed algorithm which exploits the freedom left in
\req{sigkupdate} (or \req{sigkupdate-m}).  Focusing on the
\al{OFFAR2} algorithm (that is \al{OFFAR$p$} for $p=2$), we rewrite
\req{sigkupdate} as 
\[
\sigma_k = \max[ \vartheta \nu_k, \xi_k \mu_{1,k} ]
\]
for some factor $\xi_k \in [\vartheta,1]$ which we adaptively update
as follows. We first define, for some power $\beta \in (0,1]$, an
initial gradient-norm ``target'' $t_0 =\sfrac{9}{10}\|g_k\|^\beta$ and an initial
factor $\xi_0=1$.  For $k>0$, we then update $t_k$ and
$\xi_k$, according to the rules
\beqn{rhok-rule}
\xi_k = \left\{\begin{array}{ll}
\max[\vartheta,\half \xi_{k-1}] & \tim{if } \|g_k\| \leq t_{k-1}, \\
\half(1+\xi_{k-1})& \tim{if } \|g_k\| > \max[t_{k-1},\|g_{k-1}\|] \tim{and} \xi_{k-1}< 1, \\
\xi_{k-1} & \tim{otherwise}\\
\end{array}\right.
\eeqn
and
\beqn{tk-rule}
t_k = \left\{\begin{array}{ll}
\sfrac{9}{10}\|g_k\|^\beta & \tim{if }  \|g_k\| \leq t_{k-1}, \\
t_{k-1} & \tim{otherwise.}
\end{array}\right.
\eeqn
Thus the factor $\xi_k$ decreases when the gradient's norms
converge to zero at a sufficiently fast rate (determined by the power
$\beta$), while $\xi_k$ increases to one if the gradient norms
increase\footnote{The constants $\half$ and $\sfrac{9}{10}$ may of
course be modified, but we found these to work well in our tests.}.
In addition, we choose the value of $\nu_0=\sigma_0$ with the
objective of getting $\|s_0\|$ of the order of unity\footnote{Making the
second ratio in the right-hand side of \req{skbound} equal to one.},
and set
\[
\nu_0 =\max\left[\varsigma, 6\|g_0\|\right].
\]
Finally, we use $\vartheta = 0.001$.

In what follows, we compare the standard second-order regularization algorithm
\al{AR2}\footnote{See \cite[page~65]{CartGoulToin22} with
$\eta_1=10^{-4}$, $\eta_2=0.95$, $\gamma_1 = 2 =
1/\gamma_2$,$\gamma_3=10^{20}$. $\sigma_{\min} = 10^{-4}$,
$\theta_1=0.1$ and $\theta_2=+\infty$.} with
two variants of the \al{OFFAR2} method we just described.  The first variant,
called \al{OFFAR2a} uses $\beta = 1$ in the definition of $t_0$ and
\req{tk-rule} and the second, \al{OFFAR2b}, uses $\beta = 2/3$.
All on the algorithms were run\footnote{In Matlab 
on a Dell portable computer under Ubuntu 20.04 with sixteen cores and 64
GB of memory.} on the low dimensional instances of the
problems\footnote{From their standard starting point.} of the {\sf{OPM}}
collection \cite[January 2022]{GratToin21c} listed with their dimension in
Table~\ref{testprobs}, until either $\|\nabla_x^1f(x_k)\|_2\leq \epsilon_1$, or a
maximum of 50000 iterations was reached, or evaluation of the
gradient or Hessian returned an error.

\begin{table}[htb]\footnotesize
\begin{center}
  \begin{tabular}{|l|r|l|r|l|r|l|r|l|r|l|r|}
    \hline
Problem & $n$ & Problem & $n$ & Problem & $n$ & Problem & $n$ & Problem & $n$ & Problem & $n$ \\
\hline
argauss       &  3 & chebyqad    & 10 & dixmaanl    & 12 & heart8ls   &  8 & msqrtals    & 16 & scosine     & 10 \\ 
arglina       & 10 & cliff       &  2 & dixon       & 10 & helix      &  3 & msqrtbls    & 16 & sisser      &  2 \\
arglinb       & 10 & clplatea    & 16 & dqartic     & 10 & hilbert    & 10 & morebv      & 12 & spmsqrt     & 10 \\
arglinc       & 10 & clplateb    & 16 & edensch     & 10 & himln3     &  2 & nlminsurf   & 16 & tcontact    & 49 \\
argtrig       & 10 & clustr      &  2 & eg2         & 10 & himm25     &  2 & nondquar    & 10 & trigger     &  7 \\
arwhead       & 10 & cosine      & 10 & eg2s        & 10 & himm27     &  2 & nzf1        & 13 & tridia      & 10 \\
bard          &  3 & crglvy      &  4 & eigfenals   & 12 & himm28     &  2 & osbornea    &  5 & tlminsurfx  & 16 \\
bdarwhd       & 10 & cube        &  2 & eigenbls    & 12 & himm29     &  2 & osborneb    & 11 & tnlminsurfx & 16 \\
beale         &  2 & curly10     & 10 & eigencls    & 12 & himm30     &  3 & penalty1    & 10 & vardim      & 10 \\
biggs5        &  5 & dixmaana    & 12 & engval1     & 10 & himm32     &  4 & penalty2    & 10 & vibrbeam    &  8 \\
biggs6        &  6 & dixmaanb    & 12 & engval2     &  3 & himm33     &  2 & penalty3    & 10 & watson      & 12 \\
brownden      &  4 & dixmaanc    & 12 & expfit      &  2 & hypcir     &  2 & powellbs    &  2 & wmsqrtals   & 16 \\
booth         &  2 & dixmaand    & 12 & extrosnb    & 10 & indef      & 10 & powellsg    & 12 & wmsqrtbls   & 16 \\
box3          &  3 & dixmaane    & 12 & fminsurf    & 16 & integreq   & 10 & powellsq    &  2 & woods       & 12 \\
brkmcc        &  2 & dixmaanf    & 12 & freuroth    &  4 & jensmp     &  2 & powr        & 10 & yfitu       &  3 \\
brownal       & 10 & dixmaang    & 12 & genhumps    &  5 & kowosb     &  4 & recipe      &  2 & zangwill2   &  2 \\
brownbs       &  2 & dixmaanh    & 12 & gottfr      &  2 & lminsurg   & 16 & rosenbr     & 10 & zangwill3   &  3 \\
broyden3d     & 10 & dixmaani    & 12 & gulf        &  4 & mancino    & 10 & sensors     & 10 &             &    \\
broydenbd     & 10 & dixmaanj    & 12 & hairy       &  2 & mexhat     &  2 & schmvett    &  3 &             &    \\
chandheu      & 10 & dixmaank    & 12 & heart6ls    &  6 & meyer3     &  3 & scurly10    & 10 &             &    \\
\hline
\end{tabular}
\caption{\label{testprobs} The \sf{OPM} test problems and their dimension}
\end{center}
\end{table}

Before considering the results, we recall two important comments made
in \cite{GratJeraToin22b}.  The first is that very few of the test
functions satisfy AS.3 on the whole of $\Re^n$.  While this is usually
not a problem when testing standard first-order descent methods
(because AS.3 may then be true in the level set determined by the
starting point), this is no longer the case for significantly
non-monotone methods like the ones tested here.  As a consequence, it
may (and does) happen that the gradient evaluation is attempted at a
point where its value exceeds the Matlab overflow limit, causing the
algorithm to fail on the problem.  The second comment is that the
non-monotonicity of the \al{OFFAR2} methods has another
consequence in practice: it happens on several nonconvex problems\footnote{biggs5,
biggs6, chebyqad, eg2s, hairy, indef, penalty1, penalty2, powellsq,
sensors.} that convergence of different algorithmic variants occurs to
points with gradient norm within termination tolerance (the methods
are thus achieving their objective), but these points can be far
apart and have different function values.  It is therefore
impossible to meaningfully compare the convergence performance to such
points across algorithmic variants. This reduces the set of problems
where several variants can be compared.

\begin{figure}[htb] 
\centerline{
\includegraphics[width=8cm,height=6cm]{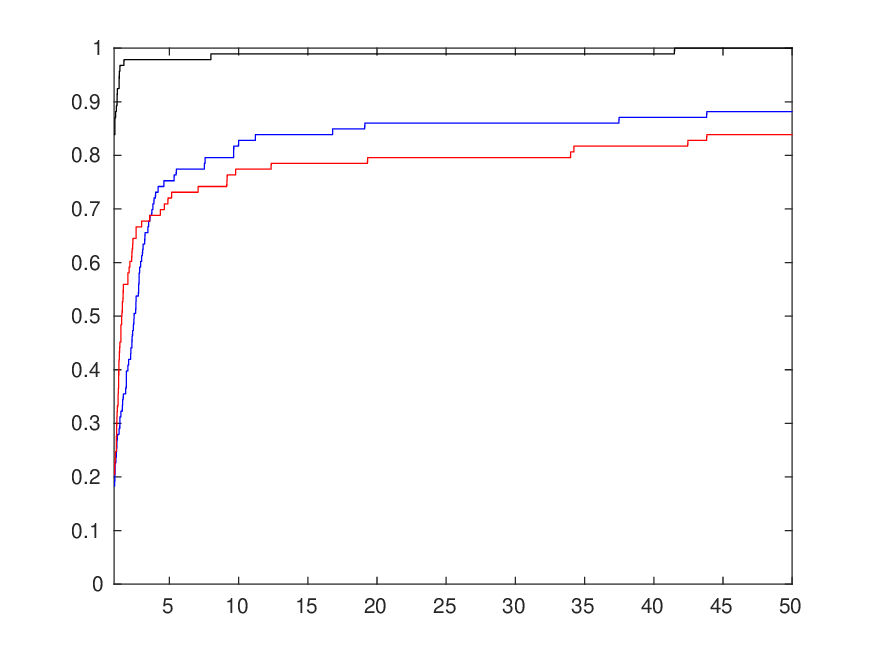}
}
\caption{\label{fig:profile}Performance profile for deterministic OFFO algorithms on
  noiseless {\sf{OPM}} problems: \al{AR2} (black), \al{OFFAR2a} (red) and \al{OFFAR2b} (blue)}
\end{figure}

We first consider the noiseless case, in which all algorithms were
terminated as soon as an iterate $x_k$ was found such that
$\|\nabla_x^1f(x_k)\|\leq 10^{-6}$.  Figure~\ref{fig:profile} shows
the corresponding performance profile \cite{DolaMoreMuns06} for the
algorithms, but we also follow \cite{PorcToin17c} and consider the
derived ``global'' measure $\pi_{\tt algo}$ to be $\sfrac{1}{50}$ of
the area below the curve corresponding to {\tt algo} in the
performance profile, for abscissas in the interval $[1,50]$. The
larger this area and closer $\pi_{\tt algo}$ to one, the closer the
curve to the right and top borders of the plot and the better the
global performance. In addition, $\rho_{\tt algo}$ denotes the
percentage of successful runs taken on all problems were comparison is
meaningful.  Table~\ref{nonoise} presents the values of the
$\pi_{\tt  algo}$ and $\rho_{\tt algo}$ statistics.  Failure\footnote{On
brownbs, cliff, cosine, curly10, himm32, genhumps, meyer3, osbornea,
recipe, scurly10, scosine, trigger, vibrbeam, yfitu.} of the
\al{OFFAR2} methods essentially occurs on ill-conditioned problems,
while \al{AR2} only fails on one problems\footnote{meyer3}.

\begin{table}[htb]
\begin{center}
\begin{tabular}{|l|r|r|r|}
  \hline
                & \al{AR2} & \al{OFFAR2a} & \al{OFFAR2b} \\
  \hline
$\pi_{\tt algo}$  &  0.99    &    0.78      &   0.83       \\
$\rho_{\tt algo}$ &  97.48   &   81.51      &   88.24      \\
\hline
\end{tabular}
\caption{\label{nonoise} Performance and reliability statistics on
  {\sf{OPM}} problems without noise}
\end{center}
\end{table}

Obviously, if we were to consider noiseless problems only, there would
be little motivation (beyond theoretical curiosity) to consider the \al{OFFAR2}
algorithms.  This is not unexpected since
\al{AR2} does exploit objective function values and we know
that identical global convergence orders may not directly translate into
similar practical behaviour. We nevertheless note that the global
behaviour of the \al{OFFAR2} algorithms is far from poor on
reasonably well-conditioned cases. 

However, as we have announced, the situation
is quite different when considering the noisy case. When noise is
present and because the value of $\mu_{1,k}$ (which is crucial
to the performance of the \al{OFFAR2} algorithnm) directly depends on
the now noisy gradient, we have modified \req{muk-def} slightly in an
attempt to attenuate the impact of noise.  We then use the
smoothed update
\[
\delta_k = \frac{9}{10}\delta_{k-1}+\frac{1}{10}\left(\frac{2\|g_k\|}{\|s_{k-1}\|^2}\right),
\ms
\mu_{1,k} = \delta_k - \theta_1\sigma_{k-1},
\]
where $\delta_0= \max[\varsigma,\|g_0\|]$, instead of
\req{muk-def}. Similarly we use a smoothed version of the gradient's
norm
\[
\tau_k  = \frac{9}{10}\tau_{k-1} +\frac{1}{10}\|g_k\|
\]
(with $\tau_{-1}= \|g_0\|$) instead of $\|g_k\|$ in \req{rhok-rule} and \req{tk-rule}.

Table~\ref{noisy} shows the reliability score $\rho_{\tt algo}$ of
the three algorithms for four increasing percentages of
relative random Gaussian noise added to the derivatives (and the
objective function for \al{AR2}), averaged over 10 runs. The tolerance
$\epsilon_1$ was set to $10^{-3}$ for these runs.

\begin{table}[htb] 
  \begin{center}
    \begin{tabular}{|l|r|r|r|r|}
      \hline
                  & $5\%$ & $15\%$ & $25\%$ & $50\%$ \\
      \hline
      \al{AR2}     & 40.67 & 30.84 & 24.54 &  6.81 \\
      \al{OFFAR2a} & 80.76 & 75.38 & 70.76 & 56.30 \\
      \al{OFFAR2b} & 85.97 & 80.67 & 72.69 & 47.98 \\
      \hline
    \end{tabular}
    \caption{\label{noisy}Reliability statistics $\rho_{\tt algo}$ for 5\%, 15\%,
      25\% and 50\% relative random Gaussian noise
      (averaged on 10 runs)}
  \end{center}
\end{table}

As anounced, the reliability of \al{AR2} decreases sharply when the
noise level increases, essentially because the decision to accept or
reject an iteration, which is based on function values, is strongly
affected by noise. By contrast, the reliability of the \al{OFFAR2}
variants decreases remarkably slowly, especially for the \al{OFFAR2a}
variant (using the more conservative $\beta=1$). That this variant is
able to solve more that 50\% of the problems with 50\% relative
Gaussian noise is quite encouraging. 

\numsection{Conclusions and Perspectives}\label{concl-s}

We have presented an adaptive regularization algorithm for nonconvex
unconstrained minimization where the
objective function is never calculated and which has, for a given
degree of used derivatives, the best-known worst-case complexity
order, not only among OFFO methods, but also among all known
optimization algorithms seeking first-order critical points.  In
particular, the algorithm using gradients and Hessians requires at
most $\calO(\epsilon_1^{-3/2})$ iterations to produce an iterate such
that $\|\nabla_x^1f(x_k)\|\leq\epsilon_1$, and at most
$\calO(\epsilon_2^{-3})$ iterations to additionally ensure that
$\lambda_{\min}[\nabla_x^2f(x_k)]\geq -\epsilon_2$.  Moreover, all
stated complexity bounds are sharp. 

These results may be extended in different ways, which
we have not included in our development to avoid too much generality
and reduce the notational burden. The first is to allow errors in
derivatives of orders 2 to $p$.  If we denote by
$\widehat{\nabla_x^if}$ the approximation of $\nabla_x^if$, it is
easily seen in the proof of Lemma~\ref{useful} that the argument
remains valid as long as, for some $\kappa_D\geq 0$,
\beqn{approxD}
\|\widehat{\nabla_x^if}(x_k)-\nabla_x^if(x_k)\| \leq \kappa_D \|s_k\|^{p+1-i}.
\eeqn
Since, for first-order analysis, the accuracy of derivatives of degree
larger than one only occurs in this lemma, we conclude that
Theorem~\ref{complexity} still hold if \req{approxD} holds.

The second extension is to replace the Euclidean norm by a more
general (possibly non-smooth) norm in $\Re^n$. As in \cite{GratToin21}, this can be achieved
without modification if first-order points are sought. When searching
for second-order points, the second-order optimality measure
$\min_{\|d\|=1}\nabla_x^2f(x_k)[d]^2$ must then be used instead of $\lambda_{\min}[\nabla_x^2f(x_k)]$.

The third extension is to weaken the gradient Lipschitz continuity
in AS.3 by only asking H\"{o}lder continuity, namely that there exist
non-negative constant $L_p$ and $\beta \in(0,1]$ such that
\beqn{LipHessian}
\|\nabla_x^p f(x) - \nabla_x^p f(y)\| \leq L_p \|x-y\|^\beta
\, \text{ for all } x,y \in \Re^n .
\eeqn
It then possible to verify that all our result remain valid with $p+1$
replaced by $p+\beta$.

Further likely generalizations include optimization in
infinite-dimensional smooth Banach spaces, a development presented for
the standard framework in \cite{GratJeraToin21}. This requires
specific techniques for computing the step and a careful handling of
the norms involved. We may also consider imposing convex constraints
on the variables \cite[Chapter~6]{CartGoulToin22}. An extension to
guarantee third-order optimality conditions (in the case where third
derivatives are available) may also be possible along the lines discussed in
\cite[Chapter~4]{CartGoulToin22}.

We also observe that the lower bound $\sigma_k \geq
\max[\mu_{1,k},\mu_{2,k}]$ has not yet, to the best of our knowledge, been used in
implementations of the standard \al{AR2} algorithm. Its incorporation in
this method is clearly possible and worth investigating.

Given the prowess of OFFO methods on noisy problems, the transition
from the present deterministic theory to the noisy context is clearly
of interest and is the object of ongoing research.  More extensive
numerical experiments involving both regularization and trust-region methods
using a variety of approximate derivatives are beyond the scope of this
paper and will be reported separately.

\appendix
\vspace*{2mm}
\flushleft{\textbf{\Large Appendix}}

\appnumsection{Proof of Theorem~\ref{sharpness2}}

\lthm{sharpness2bis}{
Let $\epsilon_2\in(0,1]$ and $p>1$.  Then there exists a $p$ times continuously
differentiable function $f_p$ from $\Re$ into $\Re$ such that the modified
\al{OFFAR$p$} applied to $f_p$ starting from the origin takes exactly
$k_\epsilon = \lceil\epsilon_2^{-\sfrac{p+1}{p-1}}\rceil$ iterations and derivative's
evaluations to produce an iterate $x_{k_\epsilon}$ such that
$|\nabla_x^1f_p(x_{k_\epsilon})|\leq \epsilon_1$
and $\lambda_{\min}[\nabla_x^2f(x_{k_\epsilon})] \geq -\epsilon_2$.
}

\proof{The proof of this result closely follows that of Theorem~\ref{sharpness}.
First select $\vartheta = 1$ (implying that $\sigma_k=\nu_k$ for all
$k$), some $\sigma_0=\nu_0>0$ and define, for all $k\in\iiz{k_\epsilon}$,
\beqn{omegak-ex-a}
\omega_k = \epsilon_2\,\frac{k_\epsilon-k}{k_\epsilon} \in [0,\epsilon_2]
\eeqn
and
\beqn{gk-ex-a}
g_k =0,
\ms
H_k = -(\epsilon_2+\omega_k)
\tim{ and }
D_{i,k} = 0,  \ms (i=3,\ldots,p),
\eeqn
so that
\beqn{grange-ex-a}
|H_k| \in[ \epsilon_2,2\epsilon_2] \subset [0,2] \tim{for all} k\in\iiz{k_\epsilon}.
\eeqn
We then set, for all $k\in\iiz{k_\epsilon}$,
\beqn{sk-ex-a}
s_k = \left(\frac{p!|H_k|}{\sigma_k}\right)^\sfrac{1}{p-1},
\eeqn
so that
\begin{align}
\sigma_k
&\eqdef \sigma_0+\sum_{j=0}^{k-1}\sigma_j|s_j|^{p+1}\label{sig-ex-a}\\
&= \sigma_0+\sum_{j=0}^{k-1}\sigma_j\left(\frac{p!|H_j|}{\sigma_j}\right)^\frac{p+1}{p-1}
= \sigma_0+(p!)^\frac{p+1}{p-1}\sum_{j=0}^{k-1}\frac{(\epsilon_2+\omega_j)^\frac{p+1}{p-1}}{\sigma_j^\frac{2}{p-1}}\nonumber\\
& \leq \sigma_0+\left(\frac{(2p!)^{p+1}}{\sigma_0^2}\right)^\frac{1}{p-1}\sum_{j=0}^{k-1}\epsilon_2^\frac{p+1}{p-1}
\leq \sigma_0+\left(\frac{(2p!)^{p+1}}{\sigma_0^2}\right)^\frac{1}{p-1}k_\epsilon\epsilon_2^\frac{p+1}{p-1}
\leq \sigma_0+2\left(\frac{(2p!)^{p+1}}{\sigma_0^2}\right)^\frac{1}{p-1}\nonumber\\
& \eqdef \sigma_{\max},\nonumber
\end{align}
where we successively used \req{sk-ex-a}, \req{gk-ex-a}, \req{omegak-ex-a} and
the definition of $k_\epsilon$.
We finally set
\[
f_0 =2^\frac{p+1}{p-1}\left(\frac{p!}{\sigma_0}\right)^\frac{2}{p-1}
\tim{and}
f_{k+1} \eqdef f_k + \half H_ks_k^2 + \sum_{i=2}^p\frac{1}{i!}D_{i,k}[s_k]^i
= f_k - \half\left(\frac{p!}{\sigma_k}\right)^\frac{2}{p-1}(\epsilon_2+\omega_k)^\frac{p+1}{p-1},
\]
yielding, using \req{sig-ex} and the definition of $k_\epsilon$, that
\[
f_0-f_{k_\epsilon}  =
\frac{1}{2}\sum_{k=0}^{k_\epsilon-1}\left(\frac{p!}{\sigma_k}\right)^\frac{2}{p-1}(\epsilon_2+\omega_k)^\frac{p+1}{p-1}
\leq 2^\frac{2}{p-1}\left(\frac{p!}{\sigma_0}\right)^\frac{2}{p-1} k_\epsilon\epsilon_2^\frac{p+1}{p}
\leq 2^\frac{p+1}{p-1}\left(\frac{p!}{\sigma_0}\right)^\frac{2}{p-1}
= f_0.
\]
As a consequence
\beqn{frange-ex-a}
f_k\in[0,f_0] \tim{for all} k\in\iiz{k_\epsilon}.
\eeqn
Observe that \req{sk-ex-a} satisfies \req{descent-m} (for the model
\req{model}), \req{gradstep-m} for $\theta_1 =1$ and \req{curvstep-m}
for $\theta_2 =1$.  Moreover \req{sig-ex-a} is the same as
\req{vkupdate-m}-\req{sigkupdate-m}. Hence the sequence $\{x_k\}$ generated by 
\[
x_0 = 0
\tim{and}
x_{k+1}=x_k+s_k
\]
may be viewed as produced by the modified \al{OFFAR$p$} algorithm given
\req{gk-ex-a}. Observe also that, for
\[
T_{k,p}(s) = f_k + g_ks + \half H_ks_k^2 +
\sum_{i=3}^p\frac{1}{i!}D_{i,k}s_k^i,
\]
one has that
\beqn{fLip-ex-a}
|f_{k+1}-T_{k,p}(s_k)|
\leq (p!)^\frac{2}{p-1}\sigma_{\max}\left(\frac{\epsilon_2+\omega_k}{\sigma_k}\right)^\frac{p+1}{p-1}
\leq \frac{\sigma_{\max}}{p!}\,|s_k|^{p+1},
\eeqn
\beqn{gLip-ex-a}
|g_{k+1}-\nabla_s^1T_{k,p}(s_k)| = |H_ks_k| = \frac{1}{p!}|s_k|^p,
\eeqn
and
\beqn{gLip-ex-a}
\begin{array}{ll}
|H_{k+1}-\nabla_s^2T_{k,p}(s_k)|
& = |H_{k+1}-H_k|
\leq |\omega_k-\omega_{k+1}|\\*[1.5ex]
& = \bigfrac{\epsilon_2}{k_\epsilon}
\leq \epsilon_2^\frac{2p}{p-1}
\leq \bigfrac{\sigma_{\max}}{\sigma_k}(\epsilon_2+\omega_k)
= \bigfrac{\sigma_{\max}}{p!}|s_k|^{p-1}
\end{array}
\eeqn
(we used $k_\epsilon \leq \epsilon_2^{-\sfrac{p+1}{p-1}}+1$ and
$\epsilon_2\leq 1$),
while, if $p>2$,
\beqn{HLip-ex-a}
|D_{i,k+1}-\nabla_s^iT_{k,p}(s_k)| = 0 \leq \frac{\sigma_{\max}}{p!}\,|s_k|^{p+1-i}
\eeqn
for $i=3,\ldots,p$.
In view of \req{grange-ex-a}, \req{frange-ex-a} and
\req{fLip-ex-a}-\req{HLip-ex-a}, we may then apply classical Hermite
interpolation to the data given by
$\{(x_k,f_k,g_k,H_k,D_{3,k},\ldots,D_{p,k})\}_{k=0}^{k_\epsilon}$ (see
\cite[Theorem~A.9.2]{CartGoulToin22} with
$
\kappa_f = \max[ 2, f_0, \sigma_{\max}/p!  ],
$
for instance) and deduce that there
exists a $p$ times continuously differentiable piecewise polynomial function $f_p$
satisfying AS.1--AS.4 and such that, for $k\in\iiz{k_\epsilon}$,
\[
f_k = f_p(x_k),
\ms
g_k = \nabla_x^1 f_p(x_k),
\ms
H_k = \nabla_x^2 f_p(x_k)
\tim{and}
D_{i,k} = \nabla_x^i f_p(x_k), \ms (i=3,\ldots,p).
\]
The sequence $\{x_k\}$ may thus be interpreted as being produced by
the \al{OFFAR$p$} algorithm applied to $f_p$ starting from $x_0=0$.
The desired conclusion then follows by observing that, from
\req{omegak-ex-a} and \req{gk-ex-a}, $g_k=0<\epsilon_1$ for all $k$ while
\[
\lambda_{\min}[H_k] = H_k < -\epsilon_2 \tim{for} k\in\iiz{k_\epsilon-1}
\tim{ and }
\lambda_{\min}[H_{k_\epsilon}]= H_{k_\epsilon} = -\epsilon_2.
\]
}

\appnumsection{Divergence of the simplified \ttal{OFFAR2} algorithm
  using \req{sigbis}}

We now prove that a modified simplified \ttal{OFFAR2} algorithm
using \req{sigbis} instead of \req{sigkupdate} may diverge. We again
proceed by constructing an example, now in $\Re^2$, where this undesirable behaviour
occur. To this aim, we first define sequences of function, gradient and Hessian
values which we will subsequently interpolate to produce the function itself.
For $k\geq 0$ and some constants $H \geq 1$ and $\theta_1\geq 1$, define
\beqn{fgHvals}
f_k = 1,
\ms
g_k = - \left(\begin{array}{c} 1 \\ 1\end{array}\right),
\ms
H_k = \left(\begin{array}{rr}
0 & 0 \\
0 & H \\
\end{array}\right)
\eeqn
and
\beqn{sigs}
\sigma_k = \frac{2(H+1)}{\sqrt{1+(H+1)^2}} < 2.
\eeqn
We may now seek the step $s_k$ which minimizes the regularized model built
from these values for given $k$, that is satisfying
\[
g_k + H_ks_k + \half \sigma_k\|s_k\| s_k = 0.
\]
Setting $[s_k]_1=1$, the first equation of this system gives that
\beqn{hsigns2}
\frac{\sigma_{k} \|s_{k}\|}{2}=\frac{\sigma_{k} \|s_{k}\|}{2} [s_k]_1 = -[g_k]_1 = 1,
\eeqn
which we may substitute in the second equation to obtain that
\[
(H+1) [s_k]_2= -[g_k]_2 = 1.
\]
This gives that
\beqn{skeven}
s_k = \left(\begin{array}{c} 1\\*[1.5ex] \bigfrac{1}{H+1} \end{array}\right)
\tim{ and thus }
\|s_k\| = \sqrt{1+ \bigfrac{1}{(H+1)^2}},
\eeqn
which is consistent with \req{sigs} and \req{hsigns2},
and we may construct a sequence $\{x_k\}$ by setting
\[
x_0 = 0
\tim{ and }
x_{k+1} = x_k + s_k
\ms (k\geq 0).
\]
Note that both  $\{[x_k]_1\}$ and $\{[x_k]_2\}$ are strictly
increasing.
We also verify, using \req{muk-def}, \req{fgHvals},
\req{hsigns2} and the bound $\theta_1\geq 1$, that
\[
\mu_{1,k}
= \frac{2\|g_k\|}{\|s_{k-1}\|^2} - \theta_1 \sigma_{k-1}
= \sqrt{2} \sigma_{k-1}  - \theta_1 \sigma_{k-1}
< \sigma_{k-1},
\]
so that \req{sigbis} gives that $\sigma_k=\sigma_{k-1}$, in accordance
with \req{sigs}.

As a consequence, the process we have just described
may be interpreted as a divergent run of the simplified \al{OFFAR2}
algo\-rithm (note that $f_k=1$ and $\|g_k\|=\sqrt{2}$ for all
$k\geq0$), provided \req{fgHvals} can be interpolated
by a function from $\Re^2$ to $\Re$ satisfying AS.1--AS.4. This is
achieved by defining $f_{k,1}=f_{k,2} = \half$ for $k \geq 0$ and
using unidimensional Hermite interpolation on the two datasets
\[
\{([x_k]_1, f_{k,1}, [g_k]_1, [H_k]_{1,1})\}_{k\geq0},
\tim{ and }
\{([x_k]_2, f_{k,2}, [g_k]_2, [H_k]_{2,2})\}_{k\geq0}.
\]
Thus we once more invoke \cite[Theorem~A.9.2]{CartGoulToin22} with
$\kappa_f = (H+1)^2$ for each dataset, which is possible because 
\[
|f_{k+1,1} -T_{k,1,2}([s_k]_1)| = |[g_k]_1[s_k]_1| = 1 = |[s_k]_1|^3 ,
\]
\[
|[g_{k+1]_1} -\nabla_s^1 T_{k,1,2}([s_k]_1)|
= |[g_{k+1}]_1 -[g_k]_1| = 0 < 1 = |[s_k]_1|^2 ,
\]
\[
|[H_{k+1}]_{1,1} -\nabla_s^2 T_{k,1,2}([s_k]_1)| = |0-0| =0 < 1  = |[s_k]_1|,
\]
and
\[
|f_{k+1,2} -T_{k,2,2}([s_k]_2) |
 = \left|-\bigfrac{1}{H+1}+\bigfrac{H}{2(H+1)^2}\right|\leq \bigfrac{1}{2(H+1)}
< (H+1)^2 |[s_k]_2|^3,
\]
\[
|[g_{k+1}]_2 -\nabla_s^1 T_{k,2,2}([s_k]_2)| =|[H_k]_{2,2}[s_k]_2| 
  = \bigfrac{H}{H+1}
< (H+1)^2|[s_k]_2|^2,
\]
\[
|[H_{k+1}]_{2,2} -\nabla_s^2 T_{k,2,2}([s_k]_2)|
=|[H_{k+1}]_{2,2} -[H_k]_{2,2} | =
0 < \frac{1}{H+1} = |[s_k]_2|,
\]
where we have defined
\[
T_{k,j,2}(s) = f_{k,j} +[g_k]_js + \half [H_k]_{j,j}[s_k]_j^2, \ms (j=1,2).
\]
We therefore deduce that there exist twice continuously differentiable
piecewise polynomial functions $f_1(x_1)$ and $f_2(x_2)$ satisfying
AS.1--AS.4 such that 
\[
f(x) = f_1(x_1)+f_2(x_2)
\]
is also twice continuously differentiable,
satisfies AS.1--AS.4 and interpolates \req{fgHvals}. This completes the argument.

\end{document}